\documentclass[a4paper,12pt]{article}
\usepackage[utf8]{inputenc}

\usepackage{amsthm,amsmath,amsfonts,amssymb}
\usepackage[english]{babel}
\usepackage{graphicx}
\usepackage{verbatim}
\usepackage{algpseudocode,algorithm,algorithmicx}
\usepackage{hyperref}
\usepackage{relsize}
\usepackage{enumitem}
\usepackage{tikz}
\usepackage{ifthen}
\usepackage{xcolor}
\usepackage[normalem]{ulem}
\usetikzlibrary{arrows, decorations.markings}
\usepackage{setspace}
\usepackage{pifont}

\definecolor{darkgreen}{RGB}{0,168,0}
\definecolor{darkblue}{RGB}{0,0,168}
\definecolor{darkred}{RGB}{168,0,0}
\definecolor{textred}{RGB}{216,0,0}

\definecolor{realred}{RGB}{160,0,0}
\definecolor{realgold}{RGB}{254,190,16}
\definecolor{realblue}{RGB}{0,0,160}

\DeclareMathOperator{\dist}{dist}

\DeclareMathOperator{\ecc}{ecc}

\DeclareMathOperator{\spec}{sp}

\newcommand{\mytilde}{\raise.17ex\hbox{$\scriptstyle\mathtt{\sim}$}}

\def\dist{\mathop{\rm dist }\nolimits}

\def\mod{\mathop{\rm mod }\nolimits}

\def\Z{\ns{Z}}

\def\Z{\ns Z}

\def\c{\mbox{\boldmath $c$}}

\def\vecu{\mbox{\boldmath $u$}}
\def\vecv{\mbox{\boldmath $v$}}

\def\vec0{\mbox{\boldmath $0$}}

\def\B{\mbox{\boldmath $B$}}

\def\B{\mbox{\boldmath $B$}}

\def\Z{\ns{Z}}

\def\Z{\mathbb Z}

\tikzset{middlearrow/.style={
		decoration={markings,
			mark= at position 0.7 with {\arrow[scale=2]{#1}} ,
		},
		postaction={decorate}
	}
}

\theoremstyle{plain}   % Cal carregar el paquet theorem.sty o amsthm.sty

\newtheorem{theorem}{Theorem}[section]
\newtheorem{proposition}[theorem]{Proposition}
\newtheorem{corollary}[theorem]{Corollary}
\newtheorem{lemma}[theorem]{Lemma}

\newtheorem{problem}[theorem]{Problem}
\theoremstyle{remark}
\newtheorem{remark}[theorem]{Remark}

\begin{document}

	\title{On bipartite $(1,1,k)$-mixed graphs
		\thanks{The research of C. Dalf\'o, M. A. Fiol and N. L\'opez has been supported by
			AGAUR from the Catalan Government under project 2021SGR00434 and MICINN from the Spanish Government under project PID2020-115442RB-I00.
			The research of M. A. Fiol was also supported by a grant from the  Universitat Polit\`ecnica de Catalunya with references AGRUPS-2022 and AGRUPS-2023.}}
	
	\author{C. Dalf\'o$^a$, G. Erskine$^b$, G. Exoo$^c$,\\ M. A. Fiol$^d$,
		J. Tuite$^b$\\
		\\
		{\small $^a$Dept. de Matem\`atica, Universitat de Lleida, Catalonia}\\
		{\small {\tt cristina.dalfo@udl.cat}},\\
		{\small $^{b}$School of Mathematics and Statistics, Open University, Milton Keynes, UK}\\
		{\small {\tt grahame.erskine,james.t.tuite@open.ac.uk}}\\
		{\small $^{c}$Dept. of Mathematics and Computer Science, Indiana State University, USA} \\
		{\small {\tt  ge@cs.indstate.edu}}\\
		{\small $^{d}$Dept. de Matem\`atiques, Universitat Polit\`ecnica de Catalunya, Barcelona, Catalonia} \\
		{\small Barcelona Graduate School of Mathematics} \\
		{\small Institut de Matem\`atiques de la UPC-BarcelonaTech (IMTech)}\\
		{\small {\tt miguel.angel.fiol@upc.edu}}}

	%\date{}

	\maketitle
	\begin{abstract}
		Mixed graphs can be seen as digraphs with arcs and edges (or digons,
		that is, two opposite arcs). In this paper, we consider the case where such
		graphs are bipartite and in which the undirected and directed
		degrees are one. The best graphs, in terms of the number of vertices, are presented for small diameters. 
		Moreover, two infinite families of such graphs with diameter $k$ and number of vertices of the order of $2^{k/2}$ are proposed, one of them being totally regular $(1,1)$-mixed graphs. In addition, we present two more infinite families called chordal ring and chordal double ring mixed graphs, which are bipartite and related to tessellations of the plane. Finally, we give an upper bound that improves the Moore bound for bipartite mixed graphs for $r = z = 1$.
	\end{abstract}
	
	
	
	%\blfootnote{
		%\begin{minipage}[l]{0.3\textwidth} \includegraphics[trim=10cm 6cm 10cm 5cm,clip,scale=0.15]{eu_logo} \end{minipage}  \hspace{-2cm} \begin{minipage}[l][1cm]{0.79\textwidth}
		%   The research of C. Dalf\'o has also received funding from the European Union's Horizon 2020 research and innovation programme under the Marie Sk\l{}odowska-Curie grant agreement No 734922.
		%  \end{minipage}}

\section{Introduction}

%It is well known that the choice of the interconnection network for a multicomputer
%or any complex system is one of the crucial problems the designer has to face.
%Indeed, the network topology largely affects the performance of the system, and it
%has an important contribution to its overall cost. As such, topologies are modeled
%by either graphs, digraphs, or mixed graphs, this has led to the following
%optimization problems:
%\begin{enumerate}
%\item
%Find graphs, digraphs, or mixed graphs of a given diameter and maximum out-degree that have a large number of vertices.
%\item
%Find graphs, digraphs, or mixed graphs of a given number of vertices and maximum out-degree that have a small diameter.
%\end{enumerate}
%
%For a more detailed description of these problems, their possible applications, the usual notation, and the theoretical background, see the comprehensive survey of Miller and \v{S}ir\'a\v{n} \cite{ms13}. For more specific results concerning mixed graphs, which are the topic of this paper, see, for example, Nguyen and Miller \cite{nm08}, and Nguyen, Miller, and Gimbert \cite{nmg07}.
%
A mixed graph can be seen as a type of digraph containing some edges (two opposite
arcs). Thus, a {\em mixed graph} $G$ with vertex set $V$ may contain (undirected) {\em
	edges} as well as directed edges (also known as {\em arcs}). From this point of
view, a {\em graph} (respectively, {\em directed graph} or {\em digraph}) has all its edges
undirected (respectively, directed). In fact, we can identify a mixed graph $G$ with its
associated digraph $G^*$ obtained by replacing all the edges with digons (two
opposite arcs or a directed $2$-cycle).
The {\em undirected degree} of a vertex $v$, denoted by $d(v)$, is the number of
edges incident to $v$. The {\em out-degree} (respectively, {\em in-degree}) of vertex $v$,
denoted by $d^+(v)$ (respectively, $d^-(v)$), is the number of arcs emanating from (respectively, to) $v$.  If $d^+(v)=d^-(v)=z$ and $d(v)=r$, for all $v \in V$, then $G$ is said to be {\em totally regular\/} of degree $(r,z)$, with $r+z=d$ (or simply {\em
	$(r,z)$-regular}).
The length of a shortest path from $u$ to $v$ is the {\it distance\/} from $u$ to $v$, and it is denoted by $\dist(u,v)$. Note that $\dist(u,v)$ may differ from $\dist(v,u)$ when the shortest paths between $u$ and $v$ involve arcs. 
%The sum of all distances from a vertex $v$, %$s(v)=\sum_{u\in V} d(v,u)$, is referred to as the {\em status\/} of $v$ (see \cite{BuckHara}). We define the {\em status vector\/} of $G$, $\mathbf{s}(G)$, as the vector constituted by the status of all its vertices. Usually, when the vector is long enough, we denote it with a short description using superscripts, that is, $\mathbf{s}(G):s_1^{n_1},s_2^{n_2},\dots,s_k^{n_k}$, where $s_1>s_2> \dots >s_k$, and $n_i$ denotes the number of vertices having $s_i$ as its local status, for all $1 \leq i \leq k$. %The {\em total status\/} of $G$, $s(G)$, is the sum of all its status:
%\[
%s(G)=\|\mathbf{s}(G)\|_1=\sum_{v\in V} s(v)=\sum_{u,v\in V} d(u,v). 
%\]
%This invariant is twice the {\em Wiener index\/} of $G$. 
The {\em out-eccentricity\/} of a vertex $u$, 
denoted by $\ecc^+(u)$, is the maximum distance from $u$ to any vertex in $G$.
Analogously, the in-eccentricity of $u$ is $\ecc^-(u)=\max\{\dist(v,u):v\in V\}$.  The maximum distance between any pair of vertices is the {\it diameter} $k=k(G)$ of $G$, that is, $k(G)=\max\{\ecc^+(u):u\in V\}$. The {\em out-radius} and {\em in-radius} of $G$ are $r^+(G)=\min\{\ecc^+(u):u\in V\}$ and $r^-(G)=\min\{\ecc^-(u):u\in V\}$, respectively.
According to Knor \cite{k96}, a {\em central vertex} is a vertex with minimum radius $r(G)=\max\{r^+(G),r^-(G)\}$. In the case of mixed graphs, we also use the concepts of {\em in-central} and {\em out-central} vertices, defined as the vertices having minimum in-radius and out-radius, respectively.
The {\em converse} of a digraph $G$, denoted by $\overline{G}$, is the digraph obtained by reversing the orientation of all arcs in $G$.
(Of course, if $G$ is a mixed graph, the edges of $\overline{G}$ remain unchanged.) Notice that  $k(G)=k(\overline{G})$, 
$r^+(G)=r^-(\overline{G})$,  and $r^-(G)=r^+(\overline{G})$.

For results concerning diameter and order (Moore bounds), see, for instance,  Nguyen and Miller \cite{nm08}, and Buset, El Amiri, Erskine, Miller, and Pérez-Rosés \cite{baemp15}.

In this paper, we consider bipartite mixed graphs with undirected and directed degrees both equal to 1. The structure of the paper is as follows. In Section~\ref{sec:11k}, we present the best (in terms of the number of vertices) bipartite $(1,1,k)$-mixed graphs for small diameters. In Section~\ref{sec:infinite-families}, we propose two infinite families of such graphs with diameter $k$ and number of vertices of the order of $2^{k/2}$, one of them being totally regular $(1,1)$-mixed graphs. In Section~\ref{sec:chordal-ring}, we give two more infinite families, called chordal ring and chordal double ring mixed graphs, which are bipartite and related to tessellations of the plane. Finally, Section~\ref{sec:upperbounds} gives an upper bound that improves the Moore bound for bipartite mixed graphs for $r = z = 1$.

\subsection{The Moore bound for bipartite $(1,1,k)$-mixed graphs}

The degree/diameter (optimization) problem for mixed graphs is the following.
\begin{problem}
	%\item[] {\bf The degree/diameter problem for mixed graphs.}
	Given three natural numbers $r,z,$
	and $k$, find the largest possible number of vertices $N(r,z,k)$ in a
	mixed graph with maximum undirected degree $r$, maximum directed out-degree $z$,
	and diameter $k$.
\end{problem}
%Buset, El Amiri, Erskine, Miller, and P\'erez-Ros\'es \cite{baemp15}
%derived a Moore-like bound of a mixed graph of diameter $k$ with
%maximum undirected degree $r$ and maximum out-degree $z$.
%An alternative approach
%for computing the bound \eqref{eq:moorebound4} has been given recently by Dalf\'o, Fiol, and L\'opez \cite{dfl16}. In order to study the bipartite case, we now use this last approach.

Here, we are interested in the case of bipartite mixed graphs. For this case, Dalf\'o, Fiol, and L\'opez \cite{dfl18} 
proved that the Moore bound for bipartite $(r,z,k)$-mixed graphs with diameter $k$ is
\begin{equation}
	\label{Moore-mix-bip}
	M_B(r,z,k)=
	2\left(A\,\frac{u_1^{k+1}-u_1}{u_1^2-1}+ B\,\frac{u_2^{k+1}-u_2}{u_2^2-1}\right),
\end{equation}
where, with $d=r+z$ and $v=(d-1)^2+4z$, 
\begin{align}
	u_1 &=\frac{d-1-\sqrt{v}}{2}, \qquad
	u_2 =\frac{d-1+\sqrt{v}}{2}, \label{u1-u2}\\
	A   &=\frac{\sqrt{v}-(d+1)}{2\sqrt{v}}, \quad \
	B   =\frac{\sqrt{v}+(d+1)}{2\sqrt{v}}. \label{A-B}
\end{align}
In the same paper \cite{dfl18}, the following results were shown:
\begin{itemize}
	\item 
	Bipartite Moore $(r,z,k)$-mixed graphs do not exist for any $r\ge 1$, $z\ge 1$,  and $k\ge 4$.
	\item
	Bipartite Moore mixed graphs with diameter $k=3$ and $r=1$ exist for any value of $z\ge  1$. In particular, some largest $(1, 1, 3)$- and $(1,1,4)$-mixed graphs were presented (see Section \ref{sec:11k}).
	\item
	There exist families of bipartite mixed graphs with diameter $k= 4, 5, 7$ and $r= 1$
	that asymptotically attain the Moore bound, for large values of $z$ being a power of a prime minus one.
\end{itemize}

%%%%%%%%%%%%%%%%%%%%%%%%%%%%%%%%%%%%%%%%%%%%%%%%%%%%%%%%%
%FIGURE (1,1,3)
%%%%%%%%%%%%%%%%%%%%%%%%%%%%%%%%%%%%%%%%%%%%%%%%%%%%%%%%%
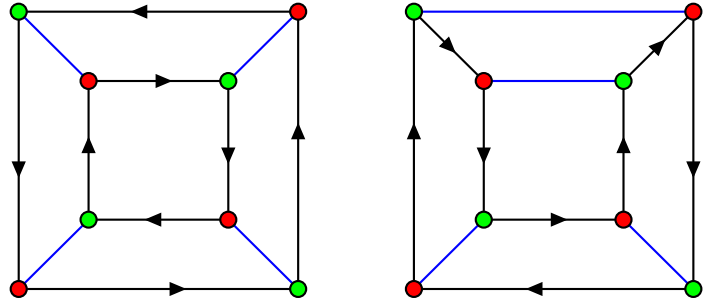
\begin{figure}[t]
	\centering
	\begin{tikzpicture}[scale=1.3]
		\tikzstyle{vertex}=[draw=black, fill=yellow!50!white, thick, shape=circle, inner sep=0, minimum height=6.0];
		\tikzstyle{edge}=[draw=blue, fill=blue, thick];
		\tikzstyle{arc}=[thick, decoration={markings,mark=at position 0.60 with {\arrow[scale=0.8,>=triangle 45]{>}}},
		postaction={decorate}];
		\newcommand\arad{2.0}
		\newcommand\brad{1.0}
		\begin{scope}[xshift=-2cm]
			\begin{scope}[rotate=45]
				\foreach \rot in {0,90,...,270}
				{
					\begin{scope}[rotate=\rot]
						\draw[arc] (0:\arad) -- (90:\arad);
						\draw[arc] (90:\brad) -- (0:\brad);
						\draw[edge] (0:\arad) -- (0:\brad);
					\end{scope}
				}
				\foreach \rot/\c in {0/red,90/green,180/red,270/green}
				{
					\begin{scope}[rotate=\rot]
						\node[vertex,fill=\c] at (0:\arad) {};
						\node[vertex,fill=\c] at (90:\brad) {};
					\end{scope}
				}
			\end{scope}
		\end{scope}
		\begin{scope}[xshift=2cm]
			\begin{scope}[rotate=45]
				\draw[arc] (90:\arad) -- (90:\brad);
				\draw[arc] (0:\brad) -- (0:\arad);
				\draw[edge] (180:\brad) -- (180:\arad);
				\draw[edge] (270:\brad) -- (270:\arad);
				\draw[edge] (0:\brad) -- (90:\brad);
				\draw[edge] (0:\arad) -- (90:\arad);
				\foreach \rot in {90,180,270}
				{
					\begin{scope}[rotate=\rot]
						\draw[arc] (90:\arad) -- (0:\arad);
						\draw[arc] (0:\brad) -- (90:\brad);
					\end{scope}
				}
				\foreach \rot/\c in {0/red,90/green,180/red,270/green}
				{
					\begin{scope}[rotate=\rot]
						\node[vertex,fill=\c] at (0:\arad) {};
						\node[vertex,fill=\c] at (90:\brad) {};
					\end{scope}
				}
			\end{scope}
		\end{scope}
	\end{tikzpicture}
	\caption{The only two bipartite Moore $(1,1,3)$-mixed graphs with $8$  vertices.}
	\label{fig1}
\end{figure}

Notice that, when $r=z=1$, \eqref{u1-u2} and \eqref{A-B} become
\begin{equation}
	u_1 =\frac{1-\sqrt{5}}{2}, \quad
	u_2 =\frac{1+\sqrt{5}}{2}, \quad
	A  =\frac{\sqrt{5}-3}{2\sqrt{5}}, \quad 
	B  =\frac{\sqrt{5}+3}{2\sqrt{5}},
	\label{u1-u2-for-rz1}
\end{equation}
and  the numbers $M_B(k)=M_B(1,1,k)$ satisfy the Fibonacci-type recurrence
$$
M_B(k)=M_B(k-1)+M_B(k-2)+2
$$
starting from $M_B(1)=2$ and $M_B(2)=4$.
In Table \ref{table1}, there are the values of $M_B(k)$ for $k=3,\ldots,16$.

%%%%%%%%%%%%%%%%%%%%%%%%%%%%%%%%%%%%%%%%%%%%%%%%%%%%%%%%%
%FIGURE (1,1,4)
%%%%%%%%%%%%%%%%%%%%%%%%%%%%%%%%%%%%%%%%%%%%%%%%%%%%%%%%%
\begin{figure}[t]
	\begin{center}
		\includegraphics[width=12cm]{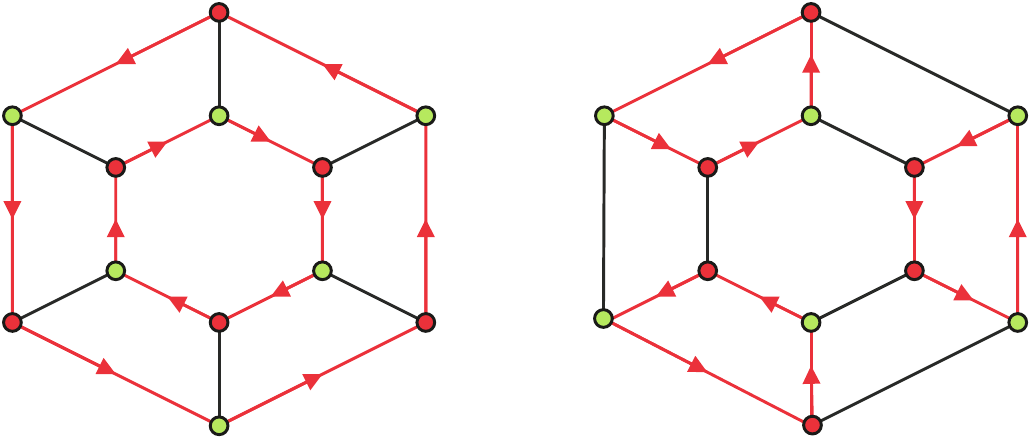}
		%\vskip -12.3cm
		\caption{Two largest bipartite $(1,1,4)$-mixed graphs (with $12$ vertices, $2$ less than the corresponding bipartite Moore bound).}
		\label{fig2}
	\end{center}
\end{figure}
%%%%%%%%%%%%%%%%%%%%%%%%%%%%%%%%%%%%%%%%%%%%%%%%%%%%%%%%%%%%%

\begin{table}[ht]
	\begin{center}
		% \begin{tabular}{|r||r|r|r|l|} \hline
			% $k$ & Lower bound & Upper bound & Moore $M_B(k)$ & $BDM(2,m)$ \\ \hline \hline
			%  3 & 8 & 8 & 8 & \\ \hline
			%  4 & 12 & 12 & 14 & \\ \hline
			%  5 & 18 & 18 & 24 & \\ \hline
			%  6 & 30 & 36 & 40 & 20 \\ \hline
			%  7 & 48  &  60 & 66 &  \\ \hline
			%  8 & 54  & 96  & 108 &  40\\ \hline
			%  9 & 96  & 158  & 176  &  \\ \hline
			% 10 & 144 & 256 & 286 &  80 \\ \hline
			% 11 & 228 & 416  & 464 &  \\ \hline
			% 12 & 312 &  674 & 742  &  160\\ \hline
			% 13 & 480 & 1092  & 1208   &  \\  \hline
			% 14 & 800 & 1766  & 1952  &  320\\ \hline
			% 15 & 1024 &  2860 & 3162  &  \\ \hline
			% 16 & 1600 &  4628 & 5116  &  640\\
			% \hline
			% \end{tabular}
		
		\begin{tabular}{|r||r|r|r|r|} \hline
			& \multicolumn{2}{|c|}{Upper bounds} & \multicolumn{2}{|c|}{Best graphs found} \\
			\hline
			$k$ & Moore $M_B(k)$ & Thm~\ref{th:james} & $BDM(2,m)$ & Computer search \\ \hline \hline
			3 & 8 & 8 & & 8* \\ \hline
			4 & 14 & 12 &  & 12* \\ \hline
			5 & 24 & 18 & & 18* \\ \hline
			6 & 40 & 36 & 20 & 30 \\ \hline
			7 & 66 &  60 &  & 48 \\ \hline
			8 & 108  & 96 & 40 & 54\\ \hline
			9 & 176 & 158  &   & 176 \\ \hline
			10 & 286 & 256 & 80 &  144 \\ \hline
			11 & 464 & 416  & & 228 \\ \hline
			12 & 742 & 674 & 160  & 312\\ \hline
			13 & 1208 & 1092  &  &  480\\  \hline
			14 & 1952 & 1766  & 320  &  800 \\ \hline
			15 & 3162 &  2860 & & 1024 \\ \hline
			16 & 5116 &  4628 & 640  &  1600\\
			\hline
		\end{tabular}

		\caption{Bounds for bipartite mixed graphs with $(r,z,k) = (1,1,k)$.}
		\label{table1}
	\end{center}
\end{table}
From \eqref{Moore-mix-bip} and \eqref{u1-u2-for-rz1}, note that the maximum possible number of vertices of a bipartite $(1,1,k)$-mixed graph is of the order $M_B(k)\sim \varphi^k$, where $\varphi$ is the golden ratio $\frac {1+\sqrt{5}}{2}\approx 1.61803$.

%%%%%%%%%%%%%%%%%%%%%%%%%%%%%%%%%%%%%%%%%%%%%%%%%%%%%%%%%
%FIGURE BDM(2,5)
%%%%%%%%%%%%%%%%%%%%%%%%%%%%%%%%%%%%%%%%%%%%%%%%%%%%%%%%%
\begin{figure}[t]
	\begin{center}
		\includegraphics[width=15cm]{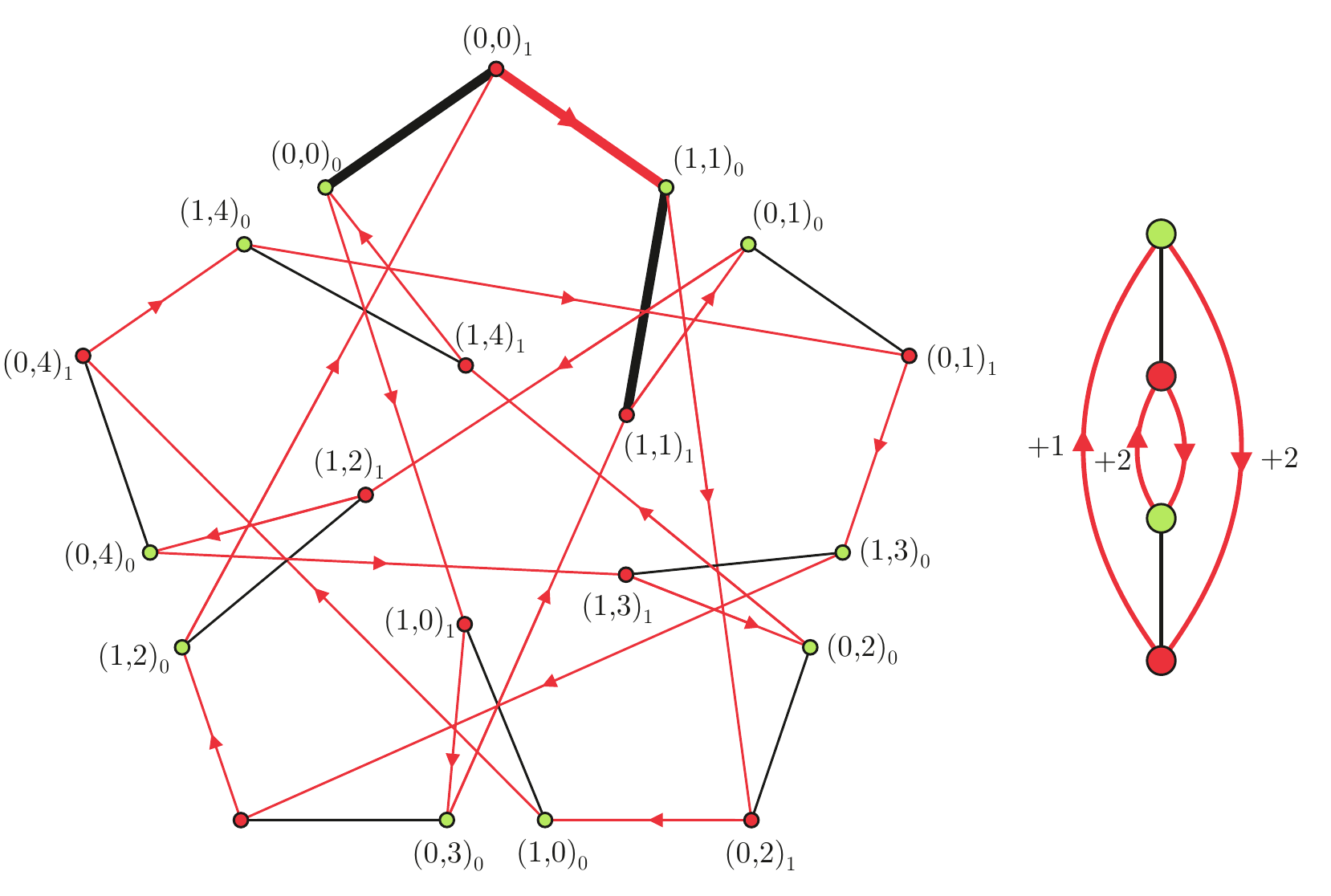}
		%\vskip -12.3cm
		\caption{The bipartite mixed graph $BDM(2,5)$ and its base graph. The thick lines in $BDM(2,5)$ represent copy $0$ of the lift.}
		\label{fig3}
	\end{center}
\end{figure}
%%%%%%%%%%%%%%%%%%%%%%%%%%%%%%%%%%%%%%%%%%%%%%%%%%%%%%%%%%%%%

% \begin{table}[t]
	%     \centering
	%     \begin{tabular}{|c||ccccccccc|}
		% \hline
		% $k$  & $2$ & $3$ & $4$ & $5$ & 6 & 7 & 8 & 9 & 10 \\
		% \hline
		% $M_B(k)$ & 4 &  8 & 14 &  24 &  40 & 66 & 108 & 176 & 286 \\
		% 				\hline
		% 			\end{tabular}
	%      \caption{The values of $M_B(k)$ for $k=2,\ldots,10$.}
	%      \label{tab:taula1}
	% \end{table}

\section{Bipartite $(1,1,k)$-mixed graphs with small diameter}
\label{sec:11k}

In this section, we present some bipartite $(1,1,k)$-mixed graphs that are best possible regarding the number of vertices. 
%for $k=3,4,5,6$ are shown in Figures \ref{fig0}, \ref{fig1} and \ref{fig2}, %respectively.
First, the mixed graphs of diameters $k=3$ and $k=4$, shown in Figures \ref{fig1} and \ref{fig2}, respectively,  were found by Dalf\'o, Fiol, and L\'opez in \cite{dfl18}. The ones with diameter $k=3$ are Moore graphs, whereas the ones with $k=4$ have order two less than the bipartite Moore bound and were proved to be the best possible. Notice that one of the mixed graphs with $k=4$ (Figure \ref{fig2}) has 2 directed 6-cycles, whereas the other has 3 directed 4-cycles.
Erskine found three other examples with the same underlying graph: one with two 6-cycles, one with 1 directed 12-cycle, and one with 1 directed 4-cycle plus 1 directed $8$-cycle. 

\begin{figure}[!ht]
	\begin{tikzpicture}[scale=0.8]
		\tikzstyle{vertex}=[draw=black, fill=yellow!50!white, thick, shape=circle, inner sep=0, minimum height=6.0];
		\tikzstyle{edge}=[draw=black, fill=blue, thick];
		\tikzstyle{arc}=[draw=darkred, ultra thick,decoration={markings,mark=at position 0.60 with {\arrow[scale=0.8,>=triangle 45]{>}}},
		postaction={decorate}];
		\newcommand\arad{4}
		\newcommand\brad{3}
		\newcommand\crad{2}
		\foreach \rot in {0,60,...,300}
		{
			\begin{scope}[rotate=\rot]
				\draw[arc] (60:\crad) -- (0:\crad);
				\draw[edge] (0:\brad) -- (0:\crad);
			\end{scope}
		}
		\foreach \rot in {0,120,240}
		{
			\begin{scope}[rotate=\rot]
				\draw[edge] (0:\arad) -- (60:\arad);
				\draw[arc] (60:\arad) -- (120:\arad);
				\draw[arc] (0:\brad) -- (60:\brad);
				\draw[arc] (0:\arad) -- (0:\brad);
				\draw[arc] (60:\brad) -- (60:\arad);
			\end{scope}
		}
		\foreach \rot/\c in {0/red,60/green,120/red,180/green,240/red,300/green}
		{
			\begin{scope}[rotate=\rot]
				\node[vertex,fill=\c] at (0:\arad) {};
				\node[vertex,fill=\c] at (60:\brad) {};
				\node[vertex,fill=\c] at (0:\crad) {};
			\end{scope}
		}
	\end{tikzpicture}
	\hskip 1cm
	%%%%%%%%%%%%%%%%%%%%%%%%%%%%%%%%%%%%%%%%%%%%%%%%%%%%%%%%%%%%%%%%
	\begin{tikzpicture}[scale=0.8]
		\tikzstyle{vertex}=[draw=black, fill=yellow!50!white, thick, shape=circle, inner sep=0, minimum height=6.0];
		\tikzstyle{edge}=[draw=black, fill=blue, thick];
		\tikzstyle{arc}=[draw=darkred, ultra thick, decoration={markings,mark=at position 0.60 with {\arrow[scale=0.8,>=triangle 45]{>}}},
		postaction={decorate}];
		\newcommand\arad{4}
		\newcommand\brad{2.5}
		\foreach \rot/\thing in {0/arc,60/arc,120/arc,180/arc,240/edge,300/arc}
		{
			\begin{scope}[rotate=\rot]
				\draw[\thing] (0:\arad) -- (60:\arad);
			\end{scope}
		}
		\foreach \rot/\thing in {0/edge,120/edge,240/arc}
		{
			\begin{scope}[rotate=\rot]
				\draw[\thing] (0:\arad) -- (15:\brad);
				\draw[\thing] (45:\brad) -- (60:\arad);
			\end{scope}
		}
		\foreach \rot in {90,210,270,330}
		{
			\begin{scope}[rotate=\rot]
				\draw[edge] (-15:\brad) -- (15:\brad);
			\end{scope}
		}
		\foreach \rot in {0,30,60,120,150,180,240,300}
		{
			\begin{scope}[rotate=\rot]
				\draw[arc] ( 15:\brad) -- (-15:\brad);
			\end{scope}
		}
		\foreach \rot in {0,120,240}
		{
			\begin{scope}[rotate=\rot]
				\draw[arc] (225:\brad) -- (75:\brad);
			\end{scope}
		}
		\foreach \rot/\c in {0/red,60/green,120/red,180/green,240/red,300/green}
		{
			\begin{scope}[rotate=\rot]
				\node[vertex,fill=\c] at (0:\arad) {};
			\end{scope}
		}
		\foreach \rot in {0,60,...,300}
		{
			\begin{scope}[rotate=\rot]
				\node[vertex,fill=green] at (15:\brad) {};
				\node[vertex,fill=red] at (45:\brad) {};
			\end{scope}
		}
	\end{tikzpicture}
	\vskip .5cm
	%%%%%%%%%%%%%%%%%%%%%%%%%%%%%%%%%%%%%%%%%%%%%%%%%%%%%%%%%%%%%%
	\begin{tikzpicture}[scale=0.8]
		\tikzstyle{vertex}=[draw=black, fill=yellow!50!white, thick, shape=circle, inner sep=0, minimum height=6.0];
		\tikzstyle{edge}=[draw=black, thick];
		\tikzstyle{arc}=[draw=darkred, ultra thick,decoration={markings,mark=at position 0.60 with {\arrow[scale=0.8,>=triangle 45]{>}}},
		postaction={decorate}];
		\newcommand\arad{4}
		\newcommand\brad{3}
		\newcommand\crad{2}
		\foreach \rot in {0,120,240}
		{
			\begin{scope}[rotate=\rot]
				\draw[edge] (0:\arad) -- (60:\arad);
				\draw[edge] (0:\crad) -- (60:\crad);
				\draw[arc] (0:\arad) -- (0:\brad);
				\draw[arc] (0:\brad) -- (0:\crad);
			\end{scope}
		}
		\foreach \rot in {60,180,300}
		{
			\begin{scope}[rotate=\rot]
				\draw[arc] (0:\arad) -- (60:\arad);
				\draw[arc] (60:\crad) -- (0:\crad);
				\draw[arc] (0:\brad) -- (0:\arad);
				\draw[arc] (0:\crad) -- (0:\brad);
			\end{scope}
		}
		\foreach \rot in {0,120,240}
		{
			\begin{scope}[rotate=\rot]
				\draw[edge] (0:\brad) to [out=225,in=45] (180:\brad);
			\end{scope}
		}
		\foreach \rot/\c in {0/red,60/green,120/red,180/green,240/red,300/green}
		{
			\begin{scope}[rotate=\rot]
				\node[vertex,fill=\c] at (0:\arad) {};
				\node[vertex,fill=\c] at (60:\brad) {};
				\node[vertex,fill=\c] at (0:\crad) {};
			\end{scope}
		}
	\end{tikzpicture}
	%%%%%%%%%%%%%%%%%%%%%%%%%%%%%%%%%%%%%%%%%%%%%%%%%%%%%%%%%%%%%%%
	\hskip 1cm
	\begin{tikzpicture}[scale=0.7]
		\tikzstyle{vertex}=[draw=black, fill=yellow!50!white, thick, shape=circle, inner sep=0, minimum height=6.0];
		\tikzstyle{edge}=[draw=black, thick];
		\tikzstyle{arc}=[ultra thick,draw=darkred,decoration={markings,mark=at position 0.60 with {\arrow[scale=0.8,>=triangle 45]{>}}},
		postaction={decorate}];
		\newcommand\arad{4}
		\foreach \rot in {0,20,...,340}
		{
			\begin{scope}[rotate=\rot]
				\draw[arc] (0:\arad) -- (20:\arad);
			\end{scope}
		}
		\foreach \rot in {0,120,240}
		{
			\begin{scope}[rotate=\rot]
				\draw[edge] (0:\arad) -- (220:\arad);
				\draw[edge] (20:\arad) -- (280:\arad);
				\draw[edge] (60:\arad) -- (200:\arad);
			\end{scope}
		}
		\foreach \rot in {0,40,...,320}
		{
			\begin{scope}[rotate=\rot]
				\node[vertex,fill=red] at (0:\arad) {};
				\node[vertex,fill=green] at (20:\arad) {};
			\end{scope}
		}
	\end{tikzpicture}
	\vskip .3cm
	%%%%%%%%%%%%%%%%%%%%%%%%%%%%%%%%%%%%%%%%%%%%%%%%%%%%%%%%
	\centering
	\begin{tikzpicture}[scale=0.7]
		\tikzstyle{vertex}=[draw=black,fill=yellow!50!white, thick, shape=circle, inner sep=0, minimum height=6.0];
		\tikzstyle{edge}=[draw=black,fill=blue,thick];
		\tikzstyle{arc}=[draw=darkred, ultra thick,decoration={markings,mark=at position 0.60 with {\arrow[scale=0.8,>=triangle 45]{>}}},
		postaction={decorate}];
		\newcommand\arad{4}
		\foreach \rot in {0,20,...,340}
		{
			\begin{scope}[rotate=\rot]
				\draw[arc] (0:\arad) -- (20:\arad);
			\end{scope}
		}
		\foreach \rot in {0,40,...,320}
		{
			\begin{scope}[rotate=\rot]
				\draw[edge] (0:\arad) -- (260:\arad);
			\end{scope}
		}
		\foreach \rot in {0,40,...,320}
		{
			\begin{scope}[rotate=\rot]
				\node[vertex,fill=red] at (0:\arad) {};
				\node[vertex,fill=green] at (20:\arad) {};
			\end{scope}
		}
	\end{tikzpicture}
	\caption{Five maximal bipartite $(1,1,5)$-mixed graph with $18$ vertices (the corresponding bipartite Moore bound is $24$).}
	\label{fig4}
\end{figure}
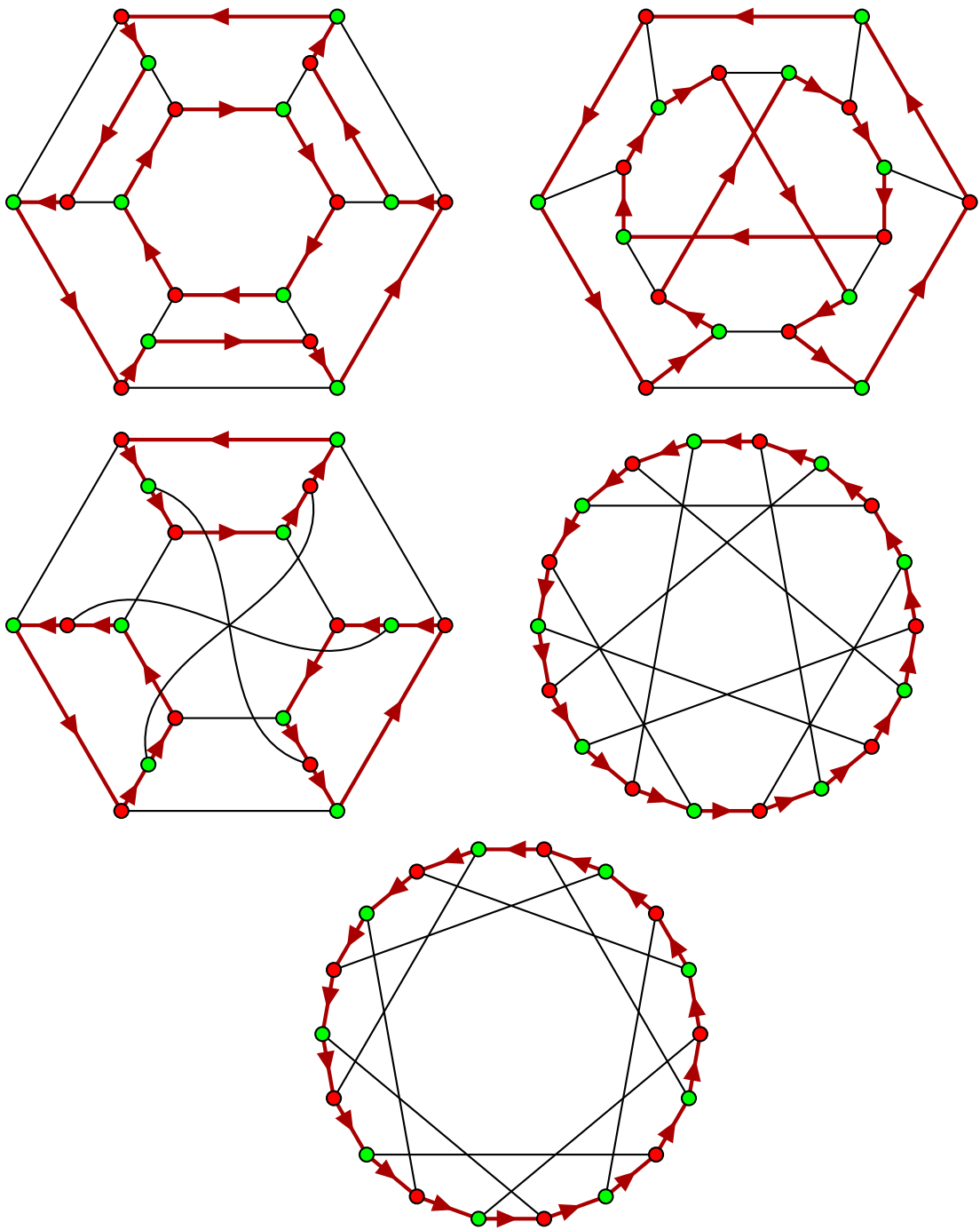

An exhaustive computational search by Exoo proved that the graphs with diameter $k=5$ and order  $N=18$ (Figure \ref{fig4}) are also the best possible. 

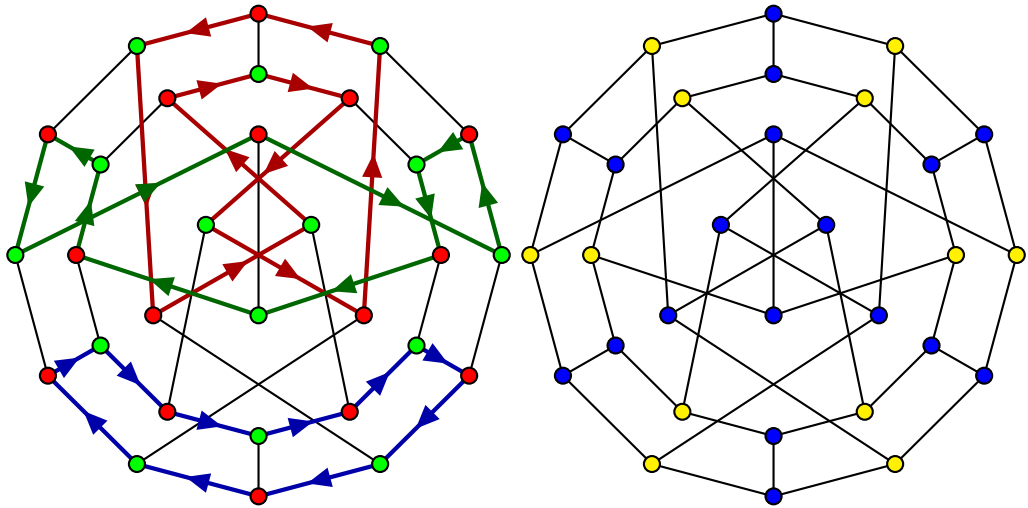
\begin{figure}[!ht]
	\begin{center}
		\begin{tikzpicture}[scale=0.8]
			\tikzstyle{vertex}=[draw=black, fill=yellow!50!white, thick, shape=circle, inner sep=0, minimum height=6.0];
			\tikzstyle{edge}=[draw=black, thick];
			\tikzstyle{arc}=[draw=darkred,ultra thick,decoration={markings,mark=at position 0.60 with {\arrow[scale=0.8,>=triangle 45]{>}}},
			postaction={decorate}];
			\newcommand\arad{4.0}
			\newcommand\brad{3.0}
			\newcommand\crad{2.0}
			\newcommand\drad{1.0}
			\begin{scope}[rotate=30]
				\foreach \rot in {0,30,...,330}
				{
					\begin{scope}[rotate=\rot]
						\draw[edge] (0:\arad) -- (30:\arad);
						\draw[edge] (0:\brad) -- (30:\brad);
					\end{scope}
				}
				\foreach \rot in {0,60,...,300}
				{
					\begin{scope}[rotate=\rot]
						\draw[edge] (0:\arad) -- (0:\brad);
					\end{scope}
				}
				\foreach \rot in {0,120,240}
				{
					\begin{scope}[rotate=\rot]
						\draw[edge] (180:\crad) -- (90:\arad);
						\draw[edge] (180:\crad) -- (270:\arad);
						\draw[edge] (0:\drad) -- (90:\brad);
						\draw[edge] (0:\drad) -- (270:\brad);
						\draw[edge] (0:\drad) -- (180:\crad);
					\end{scope}
				}
				%\uncover<2,5>{
					\draw[arc] (30:\arad) to (60:\arad);
					\draw[arc] (60:\arad) to (90:\arad);
					\draw[arc] (90:\arad) to (180:\crad);
					\draw[arc] (180:\crad) to (0:\drad);
					\draw[arc] (0:\drad) to (90:\brad);
					\draw[arc] (90:\brad) to (60:\brad);
					\draw[arc] (60:\brad) to (30:\brad);
					\draw[arc] (30:\brad) to (120:\drad);
					\draw[arc] (120:\drad) to (300:\crad);
					\draw[arc] (300:\crad) to (30:\arad);
					%}
				%\uncover<3,5>{
					\draw[arc,draw=green!40!black] (60:\crad) to (330:\arad);
					\draw[arc,draw=green!40!black] (330:\arad) to (0:\arad);
					\draw[arc,draw=green!40!black] (0:\arad) to (0:\brad);
					\draw[arc,draw=green!40!black] (0:\brad) to (330:\brad);
					\draw[arc,draw=green!40!black] (330:\brad) to (240:\drad);
					\draw[arc,draw=green!40!black] (240:\drad) to (150:\brad);
					\draw[arc,draw=green!40!black] (150:\brad) to (120:\brad);
					\draw[arc,draw=green!40!black] (120:\brad) to (120:\arad);
					\draw[arc,draw=green!40!black] (120:\arad) to (150:\arad);
					\draw[arc,draw=green!40!black] (150:\arad) to (60:\crad);
					%}
				%\uncover<4,5>{
					\draw[arc,draw=darkblue] (300:\brad) to (300:\arad);
					\draw[arc,draw=darkblue] (300:\arad) to (270:\arad);
					\draw[arc,draw=darkblue] (270:\arad) to (240:\arad);
					\draw[arc,draw=darkblue] (240:\arad) to (210:\arad);
					\draw[arc,draw=darkblue] (210:\arad) to (180:\arad);
					\draw[arc,draw=darkblue] (180:\arad) to (180:\brad);
					\draw[arc,draw=darkblue] (180:\brad) to (210:\brad);
					\draw[arc,draw=darkblue] (210:\brad) to (240:\brad);
					\draw[arc,draw=darkblue] (240:\brad) to (270:\brad);
					\draw[arc,draw=darkblue] (270:\brad) to (300:\brad);
					%}
				\foreach \rot/\c in {0,60,120,180,240,300}
				{
					\begin{scope}[rotate=\rot]
						\node[vertex,fill=red] at (0:\arad) {};
						\node[vertex,fill=green] at (0:\brad) {};
						\node[vertex,fill=green] at (30:\arad) {};
						\node[vertex,fill=red] at (30:\brad) {};
					\end{scope}
				}
				\foreach \rot/\c in {0,120,240}
				{
					\begin{scope}[rotate=\rot]
						\node[vertex,fill=red] at (180:\crad) {};
						\node[vertex,fill=green] at (0:\drad) {};
					\end{scope}
				}
			\end{scope}
		\end{tikzpicture} 
		%\vskip -1.5cm
		%\hskip .2cm
		\begin{tikzpicture}[scale=0.8]
			\tikzstyle{vertex}=[draw=black, fill=yellow!50!white, thick, shape=circle, inner sep=0, minimum height=6.0];
			\tikzstyle{edge}=[draw=black, thick];
			\tikzstyle{arc}=[decoration={markings,mark=at position 0.60 with {\arrow[scale=0.8,>=triangle 45]{>}}},
			postaction={decorate}];
			\newcommand\arad{4.0}
			\newcommand\brad{3.0}
			\newcommand\crad{2.0}
			\newcommand\drad{1.0}
			\begin{scope}[rotate=30]
				\foreach \rot in {0,30,...,330}
				{
					\begin{scope}[rotate=\rot]
						\draw[edge] (0:\arad) -- (30:\arad);
						\draw[edge] (0:\brad) -- (30:\brad);
					\end{scope}
				}
				\foreach \rot in {0,60,...,300}
				{
					\begin{scope}[rotate=\rot]
						\draw[edge] (0:\arad) -- (0:\brad);
					\end{scope}
				}
				\foreach \rot in {0,120,240}
				{
					\begin{scope}[rotate=\rot]
						\draw[edge] (180:\crad) -- (90:\arad);
						\draw[edge] (180:\crad) -- (270:\arad);
						\draw[edge] (0:\drad) -- (90:\brad);
						\draw[edge] (0:\drad) -- (270:\brad);
						\draw[edge] (0:\drad) -- (180:\crad);
					\end{scope}
				}
				\foreach \rot in {0,60,120,180,240,300}
				{
					\begin{scope}[rotate=\rot]
						\node[vertex,fill=blue] at (0:\arad) {};
						\node[vertex,fill=blue] at (0:\brad) {};
						\node[vertex,fill=yellow] at (30:\arad) {};
						\node[vertex,fill=yellow] at (30:\brad) {};
					\end{scope}
				}
				\foreach \rot in {0,120,240}
				{
					\begin{scope}[rotate=\rot]
						\node[vertex,fill=blue] at (180:\crad) {};
						\node[vertex,fill=blue] at (0:\drad) {};
					\end{scope}
				}
			\end{scope}
		\end{tikzpicture}
		\caption{A bipartite $(1,1,6)$-mixed graph with $30$ vertices (the corresponding bipartite Moore bound is $40$) and its underlying cubic graph. Blue and yellow colors indicate vertex orbits.}
		\label{fig5}
	\end{center}
\end{figure}

In Figure \ref{fig5}, we show the best mixed graph with 30 vertices and diameter $k=6$ that we found. In the figure, the 3 directed 10-cycles are shown in different colors. The question of whether or not this is the best possible solution is still open.

Table~\ref{table1} lists the largest known mixed graphs with $r=z=1$ for diameters from 3 to 16. The upper bounds are the Moore bound $M_B(k)$ from equations~(\ref{Moore-mix-bip}) and~(\ref{u1-u2-for-rz1}) and the tighter bound from Theorem~\ref{th:james} derived in Section~\ref{sec:upperbounds} below. The values shown are the orders of the graphs $BDM(2,m)$ described in Section~\ref{sec:infinite-families}, and the largest graphs found by computer search.

The computer search was exhaustive for diameters 3, 4 and 5: these values are marked with an asterisk to show that the graphs are largest possible. For diameters 7 to 16, the graphs found are a combination of Cayley graphs on groups of the stated order, and voltage lifts of a 2-vertex base graph using groups of half the stated order. Full details of the groups and generators are available from the authors on request.

%%%%%%%%%%%%%%%%%%%%%%%%%%%%%%%%%%%%%%%%%%%%%%%%%%%%%%%%%%%%%
%FIGURE (1,1,5)
%%%%%%%%%%%%%%%%%%%%%%%%%%%%%%%%%%%%%%%%%%%%%%%%%%%%%%%%%%%%%%
%%%%%%%%%%%%%%%%%%%%%%%%%%%%%%%%%%%%%%%%%%%%%%%%%%%%%%%%%%%

%%%%%%%%%%%%%%%%%%%%%%%%%%%%%%%%%%%%%%%%%%%%%%%%%%%%%%%%%
%FIGURE (1,1,6)
%%%%%%%%%%%%%%%%%%%%%%%%%%%%%%%%%%%%%%%%%%%%%%%%%%%%%%%%%

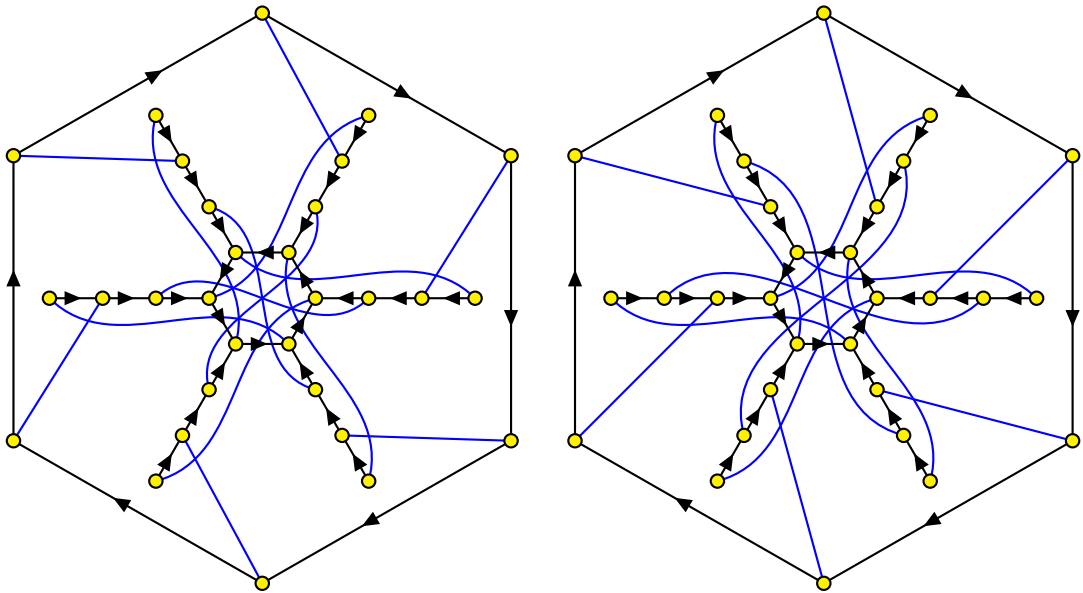
\begin{figure}[!ht]
	\centering
	\begin{tikzpicture}[scale=0.7]
		\tikzstyle{vertex}=[draw=black, fill=yellow!50!white, thick, shape=circle, inner sep=0, minimum height=5.0];
		\tikzstyle{edge}=[draw=blue, thick];
		\tikzstyle{arc}=[thick,decoration={markings,mark=at position 0.60 with {\arrow[scale=0.8,>=triangle 45]{>}}},
		postaction={decorate}];
		\newcommand\arad{4.0}
		\newcommand\brad{3.0}
		\newcommand\crad{2.0}
		\newcommand\drad{1.0}
		\newcommand\erad{5.4}
		\foreach \rot in {0,60,...,300}
		{
			\begin{scope}[rotate=\rot]
				%\uncover<3,6>{
					\draw[edge] (0:\arad) [out=135,in=315] to (120:\drad);
					%}
				%\uncover<4,6>{
					\draw[edge] (30:\erad) -- (0:\brad);
					%}
			\end{scope}
		}
		\foreach \rot in {0,60,120}
		{
			\begin{scope}[rotate=\rot]
				\draw[edge] (0:\crad) [out=225,in=45] to (180:\crad);
			\end{scope}
		}
		\foreach \rot in {0,60,...,300}
		{
			\begin{scope}[rotate=\rot]
				\draw[arc] (0:\arad) -- (0:\brad);
				\draw[arc] (0:\brad) -- (0:\crad);
				\draw[arc] (0:\crad) -- (0:\drad);
				\draw[arc] (0:\drad) -- (60:\drad);
				\draw[arc] (90:\erad) -- (30:\erad);
			\end{scope}
		}
		\foreach \rot in {0,60,...,300}
		{
			\begin{scope}[rotate=\rot]
				\node[vertex,fill=yellow] at (0:\arad) {};
				\node[vertex,fill=yellow] at (0:\brad) {};
				\node[vertex,fill=yellow] at (0:\crad) {};
				\node[vertex,fill=yellow] at (0:\drad) {};
				\node[vertex,fill=yellow] at (30:\erad) {};
			\end{scope}
		}
	\end{tikzpicture}
	\hskip .5cm
	\begin{tikzpicture}[scale=0.7]
		\tikzstyle{vertex}=[draw=black, fill=yellow!50!white, thick, shape=circle, inner sep=0, minimum height=5.0];
		\tikzstyle{edge}=[draw=blue, thick];
		\tikzstyle{arc}=[thick,decoration={markings,mark=at position 0.60 with {\arrow[scale=0.8,>=triangle 45]{>}}},
		postaction={decorate}];
		\newcommand\arad{4.0}
		\newcommand\brad{3.0}
		\newcommand\crad{2.0}
		\newcommand\drad{1.0}
		\newcommand\erad{5.4}
		\foreach \rot in {0,60,...,300}
		{
			\begin{scope}[rotate=\rot]
				\draw[edge] (0:\arad) [out=135,in=315] to (120:\drad);
				\draw[edge] (30:\erad) -- (0:\crad);
			\end{scope}
		}
		\foreach \rot in {0,60,120}
		{
			\begin{scope}[rotate=\rot]
				\draw[edge] (0:\brad) [out=225,in=45] to (180:\brad);
			\end{scope}
		}
		\foreach \rot in {0,60,...,300}
		{
			\begin{scope}[rotate=\rot]
				\draw[arc] (0:\arad) -- (0:\brad);
				\draw[arc] (0:\brad) -- (0:\crad);
				\draw[arc] (0:\crad) -- (0:\drad);
				\draw[arc] (0:\drad) -- (60:\drad);
				\draw[arc] (90:\erad) -- (30:\erad);
			\end{scope}
		}
		\foreach \rot in {0,60,...,300}
		{
			\begin{scope}[rotate=\rot]
				\node[vertex,fill=yellow] at (0:\arad) {};
				\node[vertex,fill=yellow] at (0:\brad) {};
				\node[vertex,fill=yellow] at (0:\crad) {};
				\node[vertex,fill=yellow] at (0:\drad) {};
				\node[vertex,fill=yellow] at (30:\erad) {};
			\end{scope}
		}
	\end{tikzpicture}
	\caption{Two not in-regular bipartite $(1,1,6)$-mixed graphs with $30$ vertices.}
	\label{fig6}
\end{figure}
%%%%%%%%%%%%%%%%%%%%%%%%%%%%%%%%%%%%%%%%%%%%%%%%%%%%%%%%

\section{Some infinite families of bipartite $(1,1,k)$-mixed graphs}
\label{sec:infinite-families}

For some integer $n\ge 2$, let $m=2^{n-1}+2^{n-3}$. The mixed graph $BDM(2,m)$ has independent vertex sets
\begin{align*}
	V_0=\{(0,i)_0,(1,i)_0 : i\in \Z_m\}\quad \mbox{and}\quad V_1=\{(0,i)_1,(1,i)_1 : i\in \Z_m\}.
\end{align*}
The edges are
$$
(0,i)_0\ \sim \ (0,i)_1\quad \mbox{and}\quad (1,i)_0\ \sim \ (1,i)_1,
$$
whereas the arcs are
\begin{align*}
	(0,i)_0\ \rightarrow\ (1,2i)_1\quad & \mbox{and}\quad (0,i)_1\  \rightarrow\  (1,2i+1)_0,\\
	(1,i)_0\  \rightarrow\  (0,-2i-1)_1\quad & \mbox{and}\quad (1,i)_1\  \rightarrow\  (0,-2i-2)_0,\\
\end{align*}
all with arithmetic modulo $m$. 
Thus,
$BDM(2,m)$ is a bipartite $(1,1,k)$-mixed graph on $N=4m=2^{n+1}+2^{n-1}$ vertices.
For instance, in Figures \ref{fig3} and \ref{fig8}, we show the mixed graph $BDM(2,5)$ and $BDM(2,10)$, respectively.
Notice that $BDM(2,5)$ is totally $(1,1)$-regular, but this is not the case for $BDM(2,10)$, which has vertices with in-degree $0$ and $2$. In fact, this happens for any mixed graph $BDM(2,2^{n-1}+2^{n-3})$ with $n\ge 4$ (and thus $m\ge 10$). Indeed, for $m$ even, notice that, from vertices $(0,i)_0$ and $(0,i+m/2)$, there is an arc to vertex $(1,2i)_1$. To overcome this drawback, in  Subsection \ref{sub:total-reg}, we slightly modify the adjacency rules to obtain total $(1,1)$-regularity.

In the following results, we study some of the properties of $BDM(2,m)$.

\begin{lemma}
	\label{lm:automorphism}
	Let  $\overline{0}=1$ and  $\overline{1}=0$. 
	\begin{itemize}
		\item[$(i)$] The mapping 
		$$
		\Phi_1: (\alpha,i)_{\beta}\rightarrow (\alpha,-i-1)_{\overline{\beta}},
		$$
		where $\alpha,\beta\in \Z_2$, is an involutive automorphism of $BDM(2,m)$ that interchanges its independent sets.
		\item[$(ii)$]
		The mapping 
		$$
		\Phi_2: (\alpha,i)_{\beta}\rightarrow (\alpha, i+\alpha2^{n-2}+\overline{\alpha}2^{n-3})_{\beta},
		$$
		is an automorphism of $BDM(2,m)$ of order $5$.
	\end{itemize}
\end{lemma}

\begin{proof}
	$(i)$ Let $\Gamma$ and $\Gamma^+$ denote adjacency through an edge and an arc, respectively. Since $\Phi_1$ is trivially a bijection, we only need to verify  that $\Phi_1\Gamma(\vecv)=\Gamma\Phi_1(\vecv)$ and $\Phi_1\Gamma^+(\vecv)=\Gamma^+\Phi_1(\vecv)$ for every vertex of $BDM(2,m)$:
	\begin{align*}
		\Phi_1\Gamma((\alpha,i)_{\beta})) &= \Phi_1((\alpha,i)_{\overline{\beta}}) 
		=(\alpha,-i-1)_{\beta},\\
		\Gamma\Phi_1((\alpha,i)_{\beta})) &= \Gamma((\alpha,-i-1)_{\overline{\beta}})
		=(\alpha,-i-1)_{\beta}.
	\end{align*}
	Now, assuming $\beta=0$ (the case $\beta=1$ is analogous), with $\alpha=0$,
	\begin{align*}
		\Phi_1\Gamma^+((0,i)_{0})) &=\Phi_1((1,2i)_1)=(1,-2i-1)_{0},\\
		\Gamma^+\Phi_1((0,i)_{0})) &=\Gamma^+((0,-i-1)_1)=(1,-2i-1)_{0},
	\end{align*}
	and, with $\alpha=1$,
	\begin{align*}
		\Phi_1\Gamma^+((1,i)_{0})) &=\Phi_1((0,-2i-1)_1)=(0,2i)_{0},\\
		\Gamma^+\Phi_1((1,i)_{0})) &=\Gamma^+((1,-i-1)_1)=(0,2i)_{0}.
	\end{align*}
	
	$(ii)$ We only check the directed adjacencies assuming that $\beta=1$, with $\alpha=0$,
	%so that $\Phi_2((0,i)_1)=(1,i+2^{k-3})_0$ and $\Phi_2((1,i)_1)=(1,i+2^{k-2})_1$:
	\begin{align*}
		\Phi_2\Gamma^+((0,i)_{1})) &=\Phi_2((1,2i+1)_0)=(1,2i+1+2^{n+2})_{0},\\
		\Gamma^+\Phi_2((0,i)_{1})) &=\Gamma_2^+((0,i+2^{n-3})_1)=(1,2i+2^{n-2}+1)_{0},
	\end{align*}
	and, with $\alpha=1$,
	\begin{align*}
		\Phi_2\Gamma^+((1,i)_{1})) &=\Phi_2((0,-2i-2)_1)_0)=(0,-2i-2+2^{n-3})_{0},\\
		\Gamma^+\Phi_2((1,i)_{1})) &=\Gamma^+((1,i+2^{n-2})_1)=(0,-2i-2^{n-1}-2)_{0}=(0,-2i-2+2^{n-3})_{0},
	\end{align*}
	where the last equality holds since $2^{n-3}=-2^{n-1}$ $(\mod m)$.
	
	Moreover, since, with $\alpha=0$, $\Phi_2^r((0,i)_{\beta})=(0,i+r2^{n-3})_{\beta}$, we have that $\Phi_2^r=id$ if and only if $r2^{n-3}=0$ $(\mod m)$, and $r=2^2 +1=5$ is the smallest $r$ satisfying this. The case for $\alpha=1$ is similar.
\end{proof}

\begin{figure}[!ht]
	\begin{center}
		% \vskip -1.5cm
		\includegraphics[width=8cm]{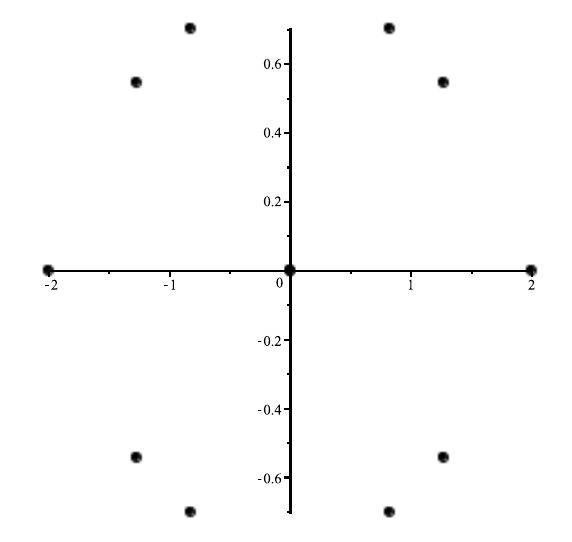}
		%\vskip -12.5cm
		\caption{The eigenvalues of the  bipartite $(1,1,6)$-mixed graph $BDM(2,5)$ in the complex plane (all of them, excepting $\pm2$, with multiplicity $2$.)}
		\label{fig7}
	\end{center}
\end{figure}

\begin{figure}[!ht]
	\begin{center}
		\includegraphics[width=14cm]{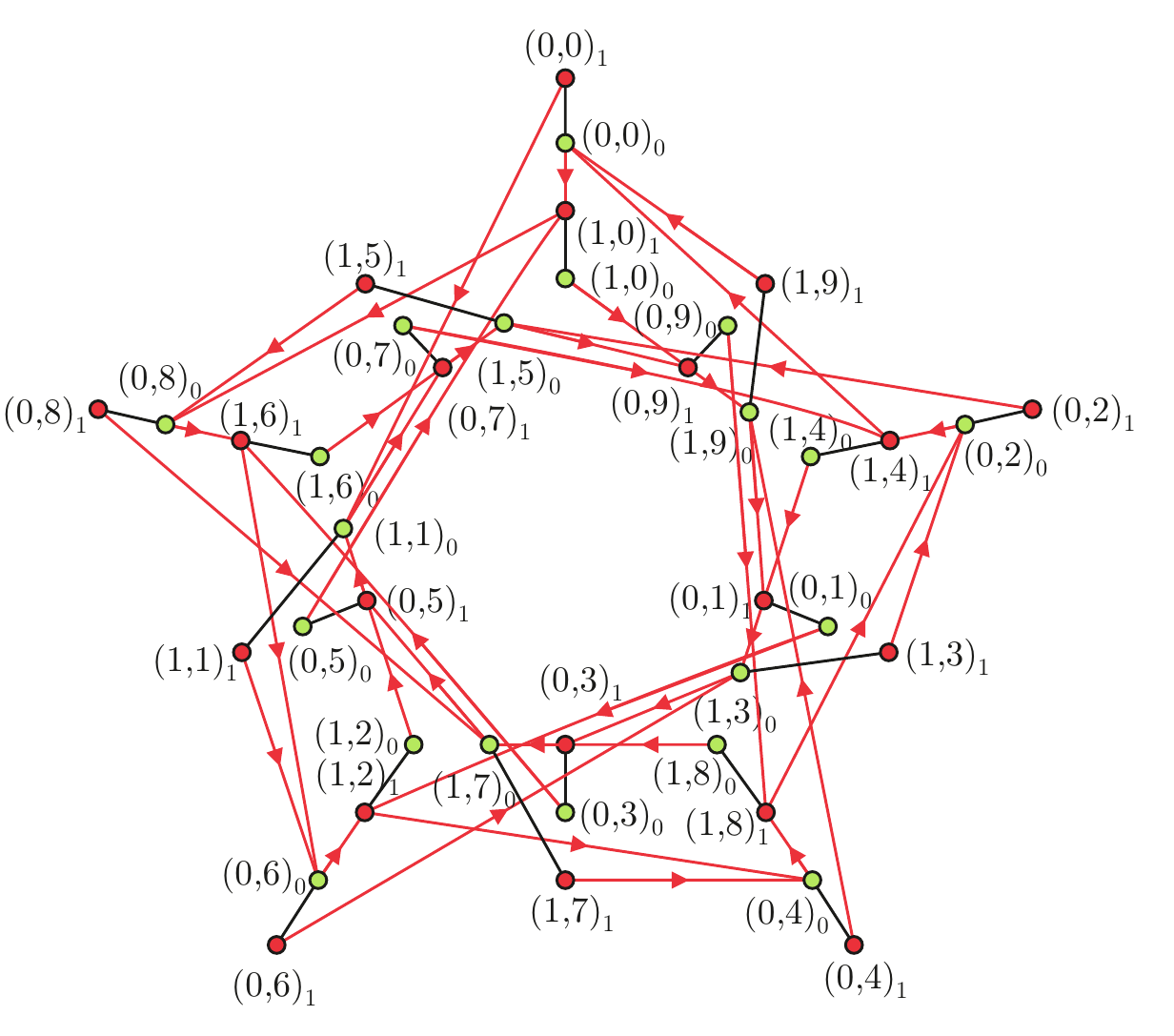}  
		%\vskip -1.5cm
		\caption{The automorphism $\Phi_2$ acting on the  bipartite $(1,1,8)$-mixed graph $BDM(2,10)$.}
		\label{fig8}
	\end{center}
\end{figure}

In Figure \ref{fig3} (left), we show the mixed graph $BDM(2,5)$ %of Figure \ref{fig4} 
drawn according to the symmetries induced by the automorphisms $\Phi_1$ and $\Phi_2$.
Notice that  $BDM(2,5)$ is totally regular. Moreover, the action of 
$\Phi_2$ allows us to construct it as a lift of the base graph with voltages on $\Z_5$ shown in Figure \ref{fig3} (right). Then, its polynomial matrix is (see Dalf\'o, Fiol,  Miller, Ryan, and \v{S}ir\'a\v{n} \cite{dfmrs19})
$$
\B=\left(
\begin{array}{cccc}
	0 & 1 & 0 & z^2\\
	1 & 0 & 1 & 0\\
	0 & z^2 & 0 & 1\\
	z & 0 & 1 & 0
\end{array}
\right),
$$
which, with $z=e^{ri\frac{2\pi}{5}}$, has the eigenvalues shown in Table \ref{table2} and Figure \ref{fig7}.

\begin{table}[t]
	\small
	\begin{center}
		\begin{tabular}{|c|cccc| }
			\hline
			$\zeta=e^{i\frac{2\pi}{5}}$, $z=\zeta^r$ & $\lambda_{r,1}$  & $\lambda_{r,2}$  & $\lambda_{r,3}$  & $\lambda_{r,4}$ \\
			\hline\hline
			$\spec(\B(\zeta^0))$ & 0 & 0 & -2 &  2 \\
			\hline
			$\spec(\B(\zeta^1))=\spec(\B(\zeta^4))$ & -0.8266-0.7015i & -0.8266+0.7015i  & 0.8266-0.7015i  & 0.8266+0.7015i   \\
			\hline
			$\spec(\B(\zeta^2))=\spec(\B(\zeta^4))$ & -1.2671-0.5445i & -1.2671+0.5445i  &  1.2671-0.5445i & 1.2671+0.5445i  \\
			\hline
		\end{tabular}
	\end{center}
	\caption{All the eigenvalues of the matrices $\B(\zeta^r)$, which yield the eigenvalues of the bipartite $(1,1,6)$-mixed graph $BDM(2,5)$.}
	\label{table2}
\end{table}

\begin{proposition}
	\label{pro:diameter}
	The diameter of the bipartite $(1,1,k)$-mixed graph $BDM(2,m)$, with $m=2^{n-1}+2^{n-3}$,  is $k=2n$.
\end{proposition}

\begin{proof}
	First, notice that when we contract all the edges of $BDM(2,m)$, we  get a bipartite  digraph with a vertex set 
	$
	V=\Z_2\times\Z_m=\{(\alpha,i): \alpha\in\Z_2, i\in \Z_m\},
	$
	and adjacencies 
	\begin{align*}
		(0,i)\ & \rightarrow\ (1,2i),\   (1,2i+1),\\
		(1,i)\ & \rightarrow\  (0,-2i-1),\  (0,-2i-2).
	\end{align*}
	This is precisely the bipartite digraph $BD(2,m)$ proposed by Fiol and Yebra in \cite{fy90}. For $m=2^{n-1}+2^{n-3}$, it was proved that $BD(2,m)$ has diameter $k=n$, see \cite[Th. 4]{fy90}.
	Thus, every path of length $\ell$ in $BD(2,m)$ induces a path of length $\ell'\le 2\ell+1$ in $BDM(2,m)$ (because every vertex of the path in  $BD(2,m)$ can turn into an edge in the corresponding path of $BDM(2,m)$).
	Then, starting from a given vertex of $BDM(2,m)$, the worst situation would be to reach a vertex through the path with adjacency pattern  $EAEAE\stackrel{(2n+1)}{\cdots\cdots} AE$ (where $E$ represents an edge, and $A$ an arc).
	From Lemma \ref{lm:automorphism}, we only need to check two initial vertices of the same partite set, say $(0,i)_1$ and $(1,i)_1$. In the first column of Tables~\ref{table3} and~\ref{table4}, we show the vertices $\vecv(n)$ and $\vecu(n)$ reached from such paths of length $2n+1$. To give a general formula, notice that the coefficient of `$i$' is clearly a power of $2$. Moreover, when we look at the adding terms of the vertices $\vecv(n)$ when $2n+1=5,9,13,17,21, \ldots$ (or the vertices $\vecu(n)$ when $2n+1=3,7,11,15,19,\ldots$)
	and omitting the signs, we have the sequence
	$1,3,13,51,205,819,\ldots$ which corresponds to A015521 in \cite{oeis} with general term $a(s)=\frac{1}{5}[4^s-(-1)^s]$ for $s=1,2,3,\ldots$ Moreover, the adding terms of the vertices $\vecv(n)$ when $2n+1=3,7,11,15,19, \ldots$ (or the vertices $\vecu(n)$ when $2n+1=1,5,9,13,17, \ldots$), omitting again the signs, are $0,2,6,26,102,\ldots$, that is, $2a(s)$ for $s=0,1,2,\ldots$. All this
	leads to the formulas shown in the tables, where the function $\phi(n)$  is given in \eqref{phi}.
	In order to show that the diameter of $BDM(2,m)$ is $k=2n$,
	in the second column of Tables \ref{table3} and \ref{table4}, there are the vertices $\vecv'(n)$ and $\vecu'(n)$ reached from the path  of the form  $AEAEAAEAE\stackrel{(2n+1)}{\cdots\cdots} AE$ and $AAEAEA\stackrel{(2n+1)}{\cdots\cdots}AE$, respectively, for $n\ge 2$. In this case, the general expression is also in the tables, with the function $\psi(n)$ given in \eqref{psi}.
	From these data, we observe that the above vertices $\vecv(n)$ and $\vecu'(n)$ can be reached, in fact, with a path of length at most $2n+1$. More precisely, when we compute the term modulo $m=2^{n-1}+2^{n-3}$, we get
	$$
	\vecv(n)=\vecv'(n-1)\quad\mbox{and}\quad \vecu(n)=\vecu'(n-3).
	$$
	Consequently, the diameter of $BDM(2,m)$ must be $k\le 2n$.
	Finally, the equality follows from the fact that, as commented, the digraph $BD(2,m)$, with $m=2^{n-1}+2^{n-3}$ has diameter $k=n$.
	
	\begin{table}[t]
		\centering
		\begin{tabular}{|c||c|c|}
			\hline
			$2n+1$  & $\vecv(n)$  & $\vecv'(n)$\\
			\hline
			3 & $(1,2i)_0$  & $(0,-4i-4)_0$ \\
			5 & $(0,-4i-1)_0$ & $(1,-8i-7)_0$\\
			7 & $(1,-8i-2)_0$ & $(0,16i+13)_0$\\
			9 &  $(0,16i+3)_0$ & $(1,32i+26)_0$\\
			11 & $(1,32i+6)_0$ & $(0,-64i-53)_0$\\ 
			13 & $(0,-64i-13)_0$ & $(1,-128i,-106)_0$ \\
			15 & $(1,-128i-26)_0$ & $(0,256i+211)_0$ \\
			17 & $(0,256i+51)_0$ & $(1,512i+422)_0$ \\
			19 & $(1,512i+102)_0$ & $\cdots$ \\
			\vdots & \vdots & \vdots \\
			$n$ even & $\left(0,\phi(n)\right)_0$ & $(1,2\phi(n))_0$\\
			$n$ odd  & $(1,2\phi(n-1))_0$ & $\left(0,\phi(n+1)\right)_0$\\
			\hline
		\end{tabular}
		\caption{The vertices $\vecv(n)$ reached from $(0,i)_1$ through a path of the form $EAEAE\stackrel{(2n+1)}{\cdots\cdots} AE$. The vertices $\vecv'(n)$ reached from $(0,i)_1$ through a path of the form $AEAEAEAE\stackrel{(2n+1)}{\cdots\cdots} AE$.}
		\label{table3}
	\end{table}
	
	\begin{table}[t]
		\centering
		\begin{tabular}{|c||c|c|}
			\hline
			$2n+1$  & $\vecu(n)$  & $\vecu'(n)$\\
			\hline
			1 & $\sim(1,i)_0$ & $\rightarrow (0,-2i-2)_0$\\
			3 & $(0,-2i-1)_0$  & $(1,-4i-4)_0$ \\
			5 & $(1,-4i-2)_0$ & $(0,-8i+7)_0$\\
			7 & $(0,8i+3)_0$ & $(1,16i+14)_0$\\
			9 &  $(1,16i+6)_0$ & $(0,-32i-29)_0$\\
			11 & $(0,-32i-13)_0$ & $(1,-64i-58)_0$\\ 
			13 & $(1,-64i-26)_0$ & $(0,128i+115)_0$ \\
			15 & $(0,128i+51)_0$ & $\cdots$ \\
			17 & $(1,256i+102)_0$ &  $\cdots$ \\
			19 & $(0,-512i-205)_0$ & $\cdots$ \\
			\vdots & \vdots & \vdots \\
			$n$ odd & $\left(0,\psi(n)\right)_0$ & $(0,\psi(n+3))_0$\\
			$n$ even  & $(1,2\psi(n-1))_0$ & $(1,2\psi(n+2))_0$\\
			\hline
		\end{tabular}
		\caption{The vertices $\vecu(n)$ reached from $(1,i)_1$ through a path of the form $EAEAE\stackrel{(2n+1)}{\cdots\cdots} AE$. The vertices $\vecu'(n)$ reached from $(1,i)_1$ through a path of the form $A$ (for $n=0$), and $AAEAEA\stackrel{(2n+1)}{\cdots\cdots}AE$ ($n\ge 1$).}
		\label{table4}
	\end{table}
	
	\vskip-.5cm
	\begin{align}
		\phi(n)&=(-1)^{\frac{n}{2}}\left[2^n i+\frac{1}{5}\left(2^n-(-1)^{\frac{n}{2}}\right)\right]\ (\mod m),%(=2^{n-1}+2^{n-3}),
		\label{phi}\\
		\psi(n)&=(-1)^{\frac{n+1}{2}}\left[2^n i+\frac{1}{5}\left(2^{n+1}-(-1)^{\frac{n+1}{2}}\right)\right]\ (\mod m).%(=2^{n-1}+2^{n-3}),
		\label{psi}
	\end{align}
	This completes the proof.
\end{proof}

%\newpage
To illustrate the situation of this proof, let us consider the case of the mixed graph $BDM(2,5)$ of Figure \ref{fig3},
the path of length $2n+1=7$ starting from $(0,i)_1$ is as follows:
\begin{align*}
	(0,i)_1\ & \sim \ (0,i)_0 \ \rightarrow \ (1,2i)_1\ \sim  \ (1,2i)_0\  \rightarrow \ (0,-4i-1)_1 \\
	& \sim\ (0,-4i-1)_0\ \rightarrow\  (1,-8i-2)_1\ \sim \  (1,-8i-2)_0,
\end{align*}
whereas the vertex $(1,-8i-2)_0$ can  be reached through the following path of length 5:
\begin{align*}
	(0,i)_1\ & \rightarrow \ (1,2i+1)_0\ \sim  \ (1,2i+1)_1\  \rightarrow \ (0,-4i-4)_0 \\
	& \sim\ (0,-4i-4)_1\ \rightarrow\  (1,-8i-7)_0=(1,-8i-2)_0.
\end{align*}
Similarly, 
the path of length $2n+1=7$ from $(1,i)_1$ is
\begin{align*}
	(1,i)_1\ & \sim \ (1,i)_0 \ \rightarrow \ (0,-2i-1)_1\ \sim  \ (0,-2i-1)_0\  \rightarrow \ (1,-4i-2)_1 \\
	& \sim\ (1,-4i-2)_0\ \rightarrow\  (0,8i+3)_1\ \sim \  (0,8i+3)_0,
\end{align*}
whereas the vertex $(0,8i+3)_0$ is, in fact, at distance $1$ from  $(1,i)_1$:
\begin{align*}
	(1,i)_1\ & \ \rightarrow \ (0,-2i-2)_0=(0,8i+3)_0.
\end{align*}

As a consequence of Proposition \ref{pro:diameter}, and in comparison with the order of the corresponding Moore bound $MD(k)\sim 1.61803^k$, we get the following result. 
\begin{corollary}
	\label{coro:order}
	The  bipartite $(1,1,k)$-mixed graph $BDM(2,m)$ has number of vertices of the order of $2^{k/2}\approx 1.4142^k$.
\end{corollary}

\subsection{Totally regular bipartite $(1,1,k)$-mixed graphs}
\label{sub:total-reg}

As commented above, the mixed graph $BDM(2,m)$ is not totally regular when $m\ge 10$. 
This subsection slightly modifies the adjacency conditions to ensure total $(1,1)$-regularity.
The bipartite mixed graph $BDM^*(2,m)$ has the same vertex set and undirected adjacencies as  $BDM(2,m)$,
whereas the arcs are now (all arithmetic modulo $m=2^{n-1}+2^{n-3}$, with $n>3$). 
\begin{itemize}
	\item
	If $i\in\left[0,m/2-1\right]$ (also as in $BDM(2,m))$:
	\begin{align*}
		(0,i)_0\ \rightarrow\ (1,2i)_1,&\quad (1,i)_0\ \rightarrow\ (0,2i+1)_1,\\
		(1,i)_0\ \rightarrow\ (0,-2i-1)_1,& \quad \ (1,i)_1\rightarrow\ (0,-2i-2)_1;
	\end{align*}
	\item
	If $i\in\left[m/2,m-1\right]$:
	\begin{align*}
		(0,i)_0\ \rightarrow\ (1,2i+1)_1,&\quad (1,i)_0\ \rightarrow\ (0,2i)_1\\
		(1,i)_0\ \rightarrow\ (0,-2i-2)_1,& \quad \ (1,i)_1\rightarrow\ (0,-2i-1)_1.
	\end{align*}
\end{itemize}

\begin{lemma}
	The bipartite graph
	$BDM^*(2,m)$ is a totally regular $(1,1,k)$-mixed graph on $N=4m=2^{n+1}+2^{n-1}$ vertices. 
\end{lemma}
\begin{proof}
	We only need to prove that every vertex $(\alpha,i)_{\beta}$ has in-degree one.
	\begin{itemize}
		\item 
		If $\alpha=\beta=0$:
		When $i$ is odd, vertex $(0,i)_0$ is adjacent from $(1,\frac{-i-1}{2})_1$;
		and when $i$ is even, vertex $(0,i)_0$ is adjacent from $(1,\frac{-i-2-m}{2})_1$.
		\item 
		If $\alpha=0$, $\beta=1$:
		When $i$ is odd, vertex $(0,i)_1$ is adjacent from $(1,\frac{-i-1-m}{2})_0$;
		and when $i$ is even, vertex $(0,i)_1$ is adjacent from $(1,\frac{-i-2}{2})_1$.
		\item 
		If $\alpha=1$, $\beta=0$:
		When $i$ is odd, vertex $(1,i)_0$ is adjacent from $(0,\frac{i-1}{2})_1$;
		and when $i$ is even, vertex $(1,i)_0$ is adjacent from $(0,\frac{i+m}{2})_1$.
		\item 
		If $\alpha=\beta=1$:
		When $i$ is odd, vertex $(1,i)_1$ is adjacent from $(0,\frac{i-1+m}{2})_0$;
		and when $i$ is even, vertex $(1,i)_1$ is adjacent from $(0,\frac{i}{2})_0$.
	\end{itemize}
\end{proof}

As in the case of $BDM(2,m)$, when we contract all the edges of $BDM^*(2,m)$, we obtain the bipartite digraph $BD(2,m)$, for $m=2^{n-1}+2^{n-3}$, with diameter $k=n$
(see again \cite{fy90}). Thus, reasoning as in the proof of Proposition \ref{pro:diameter}, the diameter of $BDM^*(2,m)$ is $k\le 2n+1$. Although, in this case, we are not able to prove that $k=2n+1$, Corollary \ref{coro:order} also applies, and the number of vertices of $BDM^*(2,m)$ is also of the order of $\sqrt{2}^k$.

\section{Chordal ring mixed graphs}
\label{sec:chordal-ring}

In this section, we present a family of bipartite $(1,1,k)$-mixed graphs that are related to tessellations of the plane  (see Yebra, Fiol, Morillo, and Alegre \cite{yfma85}).
Let $n\ge 2$ and $c<n$  be, respectively, even and odd numbers. The chordal ring mixed graph $CRM(n,c)$ is a mixed graph with vertex set $V=\Z_n$ (all arithmetic will be modulo $n$), with arcs $i\rightarrow i+1$ (forming a directed cycle) and edges $i\sim i+c$ if $i$ is odd (these are the `chords').
Given the diameter $k$, we want to find the value of $c$ such that the graph $CRM(n,c)$ has the maximum number of vertices.
Arden and Lee studied this problem \cite{al81}, and Yebra, Fiol,  Morillo, and Alegre \cite{yfma85} in the case of undirected graphs, which were called `chordal ring networks'. In fact, such graphs were already studied in another context by Coxeter \cite{c50}.

In our case of $(1,1,k)$-mixed graphs, the following result gives a Moore-like bound for their number of vertices.
\begin{lemma}
	The maximum number of vertices of a $CRM(n,c)$ of a bipartite chordal ring mixed graph with diameter $k$ is
	\begin{equation}
		\left\{
		\begin{array}{cc}
			\frac{1}{2}(k+1)^2 & \mbox{if $k$ is odd,}\\[.2cm]
			\frac{1}{2}k(k+2) & \mbox{if $k$ is even.}
		\end{array}    
		\right.
		\label{CRM(k)}
	\end{equation}
\end{lemma}
\begin{proof}
	From the adjacency conditions, we observe that there are at most $d+1$ vertices at distance $d\ge 0$ from vertex 0 (that is, $\{0\}$, $\{-c,1\}$, $\{-c+1,2,1+c\}$, $\{-2c+2,-c+2,3,2+c\}$, \ldots ). Then, if the mixed graph is bipartite with odd diameter, the maximum number of vertices is bounded above by twice the number of even vertices of a
	$CRM(n,c)$, which is equal to $2(1+3+\cdots+k)=\frac{1}{2}(k+1)^2$.
	Similarly, if the diameter is even, we have 
	$2(2+4+\cdots+k=\frac{1}{2}k(k+2)$.
\end{proof}

Next, we show that the upper bound in \eqref{CRM(k)} can be attained if $k$ is odd,
but not when $k$ is even.

% In  the following result we show that the upper bound $CRM(k)$ is attained when the diameter $k$ is odd. Moreover, we construct some families with number of vertices of the order of $O(k^2)$ when the diameter $k$ is even.

\begin{theorem}
	\label{th:CRM}
	\begin{itemize}
		\item[$(a)$]
		If $k$ is odd, $k=2\ell-1$ with $\ell\ge 2$, there exists a chordal ring mixed graph $CRM(n,c)$ with diameter $k$, order 
		$n=2\ell^2=\frac{1}{2}(k+1)^2$, and chordal length $c=2\ell-1=k$.
		\item[$(b)$]
		If $k$ is even, with $k\equiv 0$ $(\mod 4)$, $k=2\ell=4t$ with $t\ge 1$, there exists a chordal ring mixed graph $CRM(n,c)$ with diameter $k$, order 
		$n=8t^2+2=\frac{1}{2}k^2+2$, and chordal length $c=(2t-1)^2+2t=(\frac{k}{2}-1)^2+\frac{k}{2}$.
		\item[$(c1)$]
		If $k$ is even, $k\equiv 6$ $(\mod 8)$, $k=2\ell=8t-2$ with $t\ge 1$, there exists a chordal ring mixed graph $CRM(n,c)$ with diameter $k$, order 
		$n=2(4t-3)(4t-4)+4=k(\frac{k}{2}-1)+4$, and chordal length $c=8t^2-8t+3=\frac{1}{8}(k+2)^2-k+1$.
		\item[$(c2)$]
		If $k$ is even, $k\equiv 2$ $(\mod 8)$, $k=2\ell=8t-6$ with $t\ge 1$, there exists a chordal ring mixed graph $CRM(n,c)$ with diameter $k$, order 
		$n=2(4t-3)(4t-4)+4=k(\frac{k}{2}-1)+4$, and chordal length $c=24t^2-44t+23=\frac{3}{8}(k+6)^2-\frac{11}{2}k-10$.
	\end{itemize}
\end{theorem}
% Given any given odd  diameter $k=2\ell-1$, with $\ell\ge 2$,  the mixed graph $CRM(n,c)$ has maximum possible order $n=2\ell^2=\frac{1}{2}(k+1)^2$ when its chord length is  $c=2\ell-1$.

\begin{proof}
	Consider the plane divided into unit squares and divide each square by its anti-diagonal line, forming two right triangles.
	Associate to each odd vertex $i$  and even vertex $i+c$ (forming an edge) an upper and its adjacent lower triangle. The vertices $i\pm 1$ adjacent to and from $i$ are represented with the adjacent triangles, as shown in Figure \ref{fig10}. Following this procedure, the vertices can be arranged in a planar pattern, as shown in the same figure. Notice that,
	when every triangle of the regular tessellation of the plane receives a number modulo $n$ according to the adjacency rules, 
	the distribution of these numbers in the plane
	repeats itself periodically. This fact is illustrated again in Figure \ref{fig10} for the graph $CRM(32,7)$. Moreover, $n$ triangles form every tile, and they {\bf periodically tessellate the plane.}
	Stated in this context, our concern is to find and `construct' a given  tile (that is, find an integer $c$ that generates it) that tessellates the plane and has a maximum area (or number of unit triangles) for a
	given diameter $k$.\\
	A more precise approach considers that the automorphism group of $CRM(n,c)$ has two orbits constituted by the even and odd vertices. Consequently, we must compute the number of vertices from, say, $0$ and $-c$. This is shown in Figure \ref{fig9} and
	Table \ref{table5}, where we indicate the vertices whose maximum distance to $0$ or $-c$ is $1,2,\ldots,8$. If we sum up odd diameters, we obtain the bounds in \eqref{CRM(k)} again.
	% but for diameter $k$ even, we get the bound
	% $CRM(k)= 2+2(4+6+8+\cdots+k)+(k+2)$ if $k$ is even.
	
	\begin{figure}[!ht]
		\begin{center}
			\includegraphics[width=9cm]{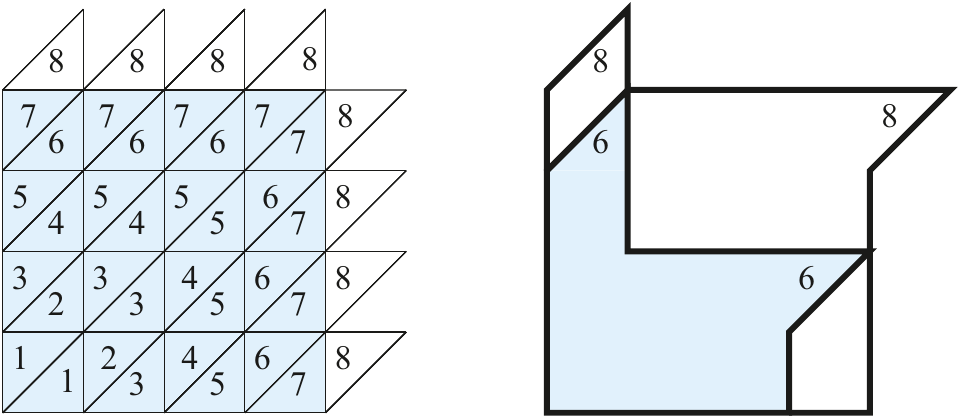}  
			%\vskip -1.5cm
			\caption{Vertices at a given distance $1,2,\ldots,8$ from $0$ or $-c$, and optimal tiles.}
			\label{fig9}
		\end{center}
	\end{figure}
	
	\begin{table}[ht]
		$$
		\begin{array}{|c|cccccccc|}
			\hline
			\max-\dist  &   &    &      &  {\rm vertices}   &      &      &   &   \\ 
			\hline
			1 &   &     &       &   -c+1   &   1 &      &       &      \\
			2 &   &     &       & -c+2   &   2  &     &      &      \\
			3 &   &     &   -2c+2    & -c+3 &   3  &  2+c &      &      \\
			4 &   &    & -2c+3 & -c+4 &   4  &  3+c &      &      \\
			5 &  &  -3c+3    & -2c+4 & -c+5 &   5  &  4+c & 3+2c &      \\
			6 &  & -3c+4 & -2c+5 & -c+6 &   6  &  5+c & 4+2c &     \\
			7 &-4c+4  & -3c+5 & -2c+6 & -c+7 &   7  &  6+c & 5+2c & 4+3c \\
			8 & -4c+5 & -3c+6 & -2c+7 & -c+8& 8& 7+c& 6+2c& 5+3c\\
			\hline
		\end{array}
		$$
		\caption{The vertices $i$ of $CRM(n,c)$ such that $\max\{\dist(0,i),\dist(-c,i)\}$ is $\max-\dist=1,2,\ldots,8$.}
		\label{table5}
	\end{table}
	
	$(a)$ To show that for odd diameter $k$ the bound on the order is attained, 
	the appropriate tile is a square of length $\ell$, as it is shown for $k=7$ in Figures \ref{fig9} (left shaded area) and Figure \ref{fig10}, with $n=2\ell^2$ triangles (vertices) and diameter $k=2\ell-1$, together with its tessellation.
	% In the bottom right square, there are the vertices at successive distance $1,2,\ldots,7$ from the even vertex $0$ (by the symmetry of the tile, 
	% the distances from an odd vertex are obtained by reflection with respect to the secondary diagonal). 
	
	\begin{figure}[!ht]
		\begin{center}
			\includegraphics[width=9cm]{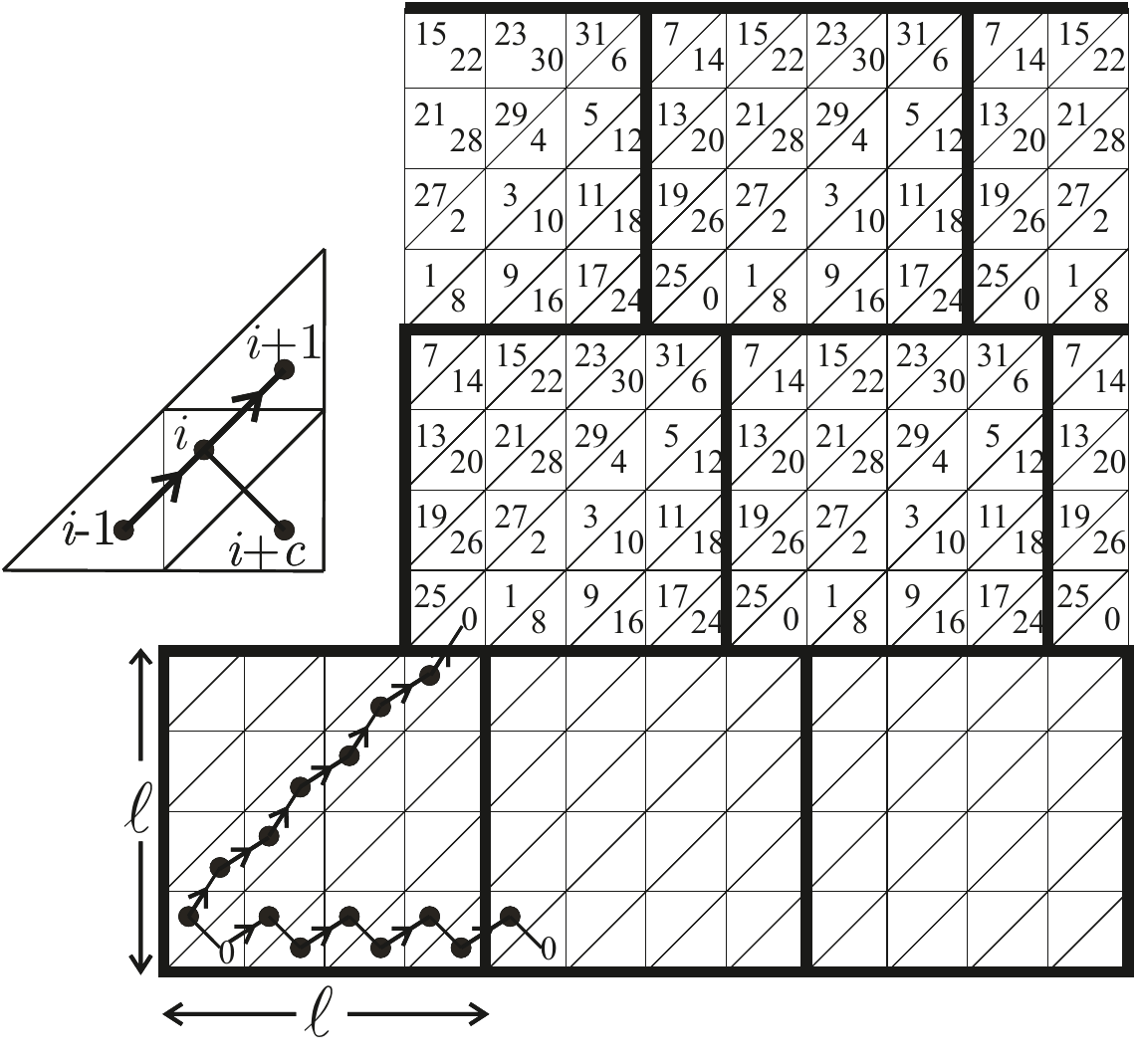}  
			%\vskip -1.5cm
			\caption{The plane tessellation corresponding to the chordal ring mixed graph  $CRM(32,7)$, with diameter $k=7$.}
			\label{fig10}
		\end{center}
	\end{figure}
	
	It remains to show that a suitable choice of $c$ can generate such a tile. For this, note that 
	these values produce the given periodic pattern, which is characterized by the position of the `zeros'. 
	To obtain this distribution, we have to express the null effect of translations along two linearly independent vectors
	(each of them associated with a path as shown in Figure \ref{fig10}) that generate the pattern. Choosing them as in the figure, $c$ must satisfy
	\begin{align*}
		\ell+\ell c &\equiv0 \ (\mod 2\ell^2) \quad\mbox{(horizontal path)},\\
		2\ell-1-c &\equiv 0 \ (\mod 2\ell^2)\quad \mbox{(diagonal path)},
	\end{align*}
	with trivial solution $c=2\ell-1$, as claimed.\\
	
	$(b)$ When the diameter is even, of the form $k=2\ell=4t$, see, for instance, Figure \ref{fig9} (left) for $k=8$, the optimal tiles clearly do not tessellate because of the bordering triangles. The best we can do is to remove all such triangles except two. This results in the tile of Figure \ref{fig9} (on the right), with order $n=2\ell^2+2$, or its equivalent (because of periodicity) L-shaped tile of Figure \ref{fig11}.
	
	\begin{figure}[!ht]
		\begin{center}
			\includegraphics[width=7cm]{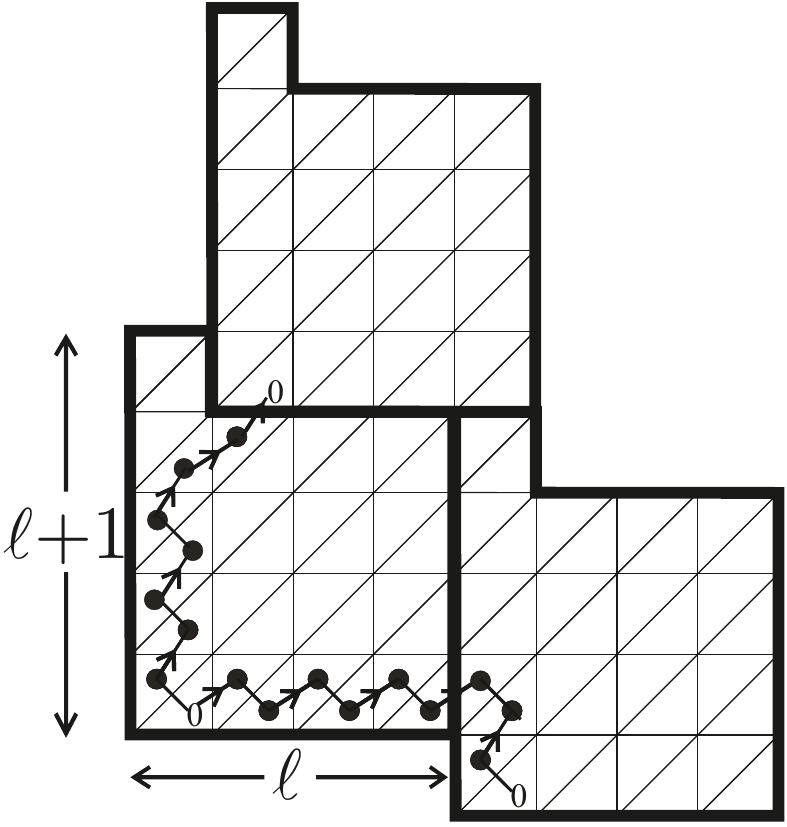}  %\vskip -1.5cm
			\caption{The plane tessellation corresponding to the chordal ring mixed graph  $CRM(34,13)$, with diameter $k=8$.}
			\label{fig11}
		\end{center}
	\end{figure}
	
	Then, the equations to obtain the right value of $c$ are
	\begin{align*}
		(\ell-1)+(\ell+1)c &\equiv 0 \ (\mod n),\\
		(\ell+1)-(\ell-1)c&\equiv 0 \ (\mod n),
	\end{align*}
	or, in matrix form,
	$$
	\left(
	\begin{array}{cc}
		\ell-1 & \ell+1\\
		\ell+1 & -\ell+1
	\end{array}
	\right)
	\left(
	\begin{array}{c}
		1\\
		c
	\end{array}
	\right)=
	n\left(
	\begin{array}{c}
		\alpha\\
		\beta
	\end{array}
	\right).
	$$
	Solving the system 
	(notice that the determinant of the $2\times 2$ matrix is 0 $(\mod n)$), we have
	\begin{align*}
		1  & = (\ell-1)\alpha+(\ell+1)\beta=(2t-1)\alpha+(2t+1)c,\\
		c &=(\ell+1)\alpha-(\ell-1)\beta= (2t+1)\alpha-(2t-1)\beta.
	\end{align*}
	Since $c$ must be an (odd) integer, a solution is obtained by taking $\alpha=t$ and $\beta=-t+1$, so that the first equation holds and $c=(2t-1)^2+2t$, as claimed.
	
	$(c)$ When the diameter is even, of the form $k=2\ell=8t-2$ (case $(c1)$) or $k=2\ell=8t-6$ (case $(c2)$), we could use, in principle, the same tile as in $(b)$. However, the corresponding tessellation yields no solution with `step' 1, and, hence, we do not obtain a proper chordal ring mixed graph (see the remark after this proof). Then, the best solution is the right shaded tile shown in Figure~\ref{fig9} for $k=6$ (or the corresponding tile for $k=10$ in Figure \ref{fig12}) with order $n=2\ell(\ell-1)+4$. 
	
	\begin{figure}[!ht]
		\begin{center}
			\includegraphics[width=8cm]{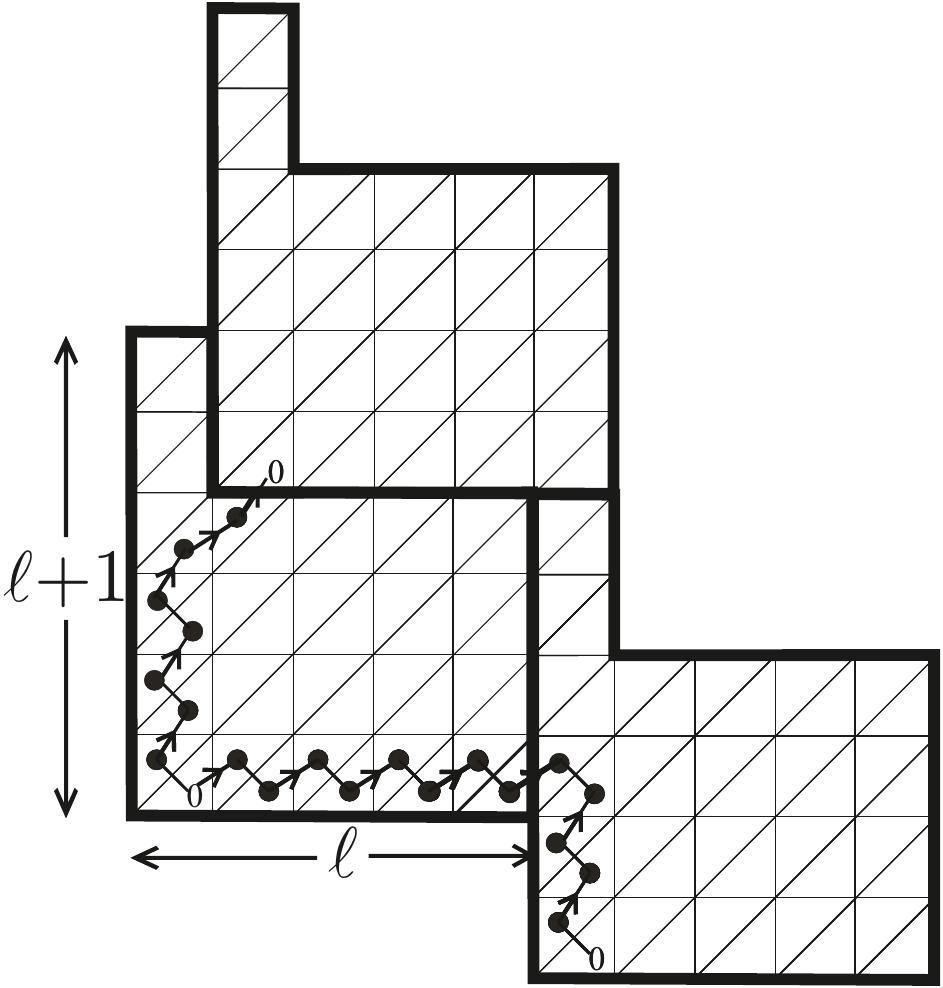}  
			%\vskip -1.5cm
			\caption{The plane tessellation corresponding to the chordal ring mixed graph  $CRM(44,31)$, with diameter $k=10$.}
			\label{fig12}
		\end{center}% \end{figure}
\end{figure}

Then, in both cases, the equations to obtain the right value of $c$ are 
$$
\left(
\begin{array}{cc}
	\ell-2 & \ell+2\\
	\ell & -\ell+2
\end{array}
\right)
\left(
\begin{array}{c}
	1\\
	c
\end{array}
\right)=
n\left(
\begin{array}{c}
	\alpha\\
	\beta
\end{array}
\right)\quad \Rightarrow\quad 
\left\{\begin{array}{ll}
	1 \!\!\!&= (\ell-2)\alpha+(\ell+2)\beta,\\
	c\!\!\! &=\ell\alpha+(-\ell+2)\beta.
\end{array}
\right.
$$
In  case $(c1)$, we get $\alpha=t$ and $\beta=-t+1$; whereas in the case $(c2)$, we have $\alpha=3t-1$ and $\beta=-3t+4$.
Then, depending on the values $\ell=\frac{1}{2}(8t-2)$ or $\ell=\frac{1}{2}(8t-6)$, we obtain the solutions $c=8t^2-87+3$ and $c=24t^2-44+23$.
\end{proof}

% \begin{figure}[!ht]
%  	\begin{center}
	% 		\includegraphics[width=6cm]{tessellation(c).pdf}  
	%  %\vskip -1.5cm
	%  	\caption{The plane tessellation corresponding to the mixed chordal ring graph  $CRM(44,31)$, with diameter $k=10$.}
	%  		\label{fig13}
	%  	\end{center}% \end{figure}
	%   \end{figure}

\begin{figure}[!ht]
\begin{center}
\includegraphics[width=14cm]{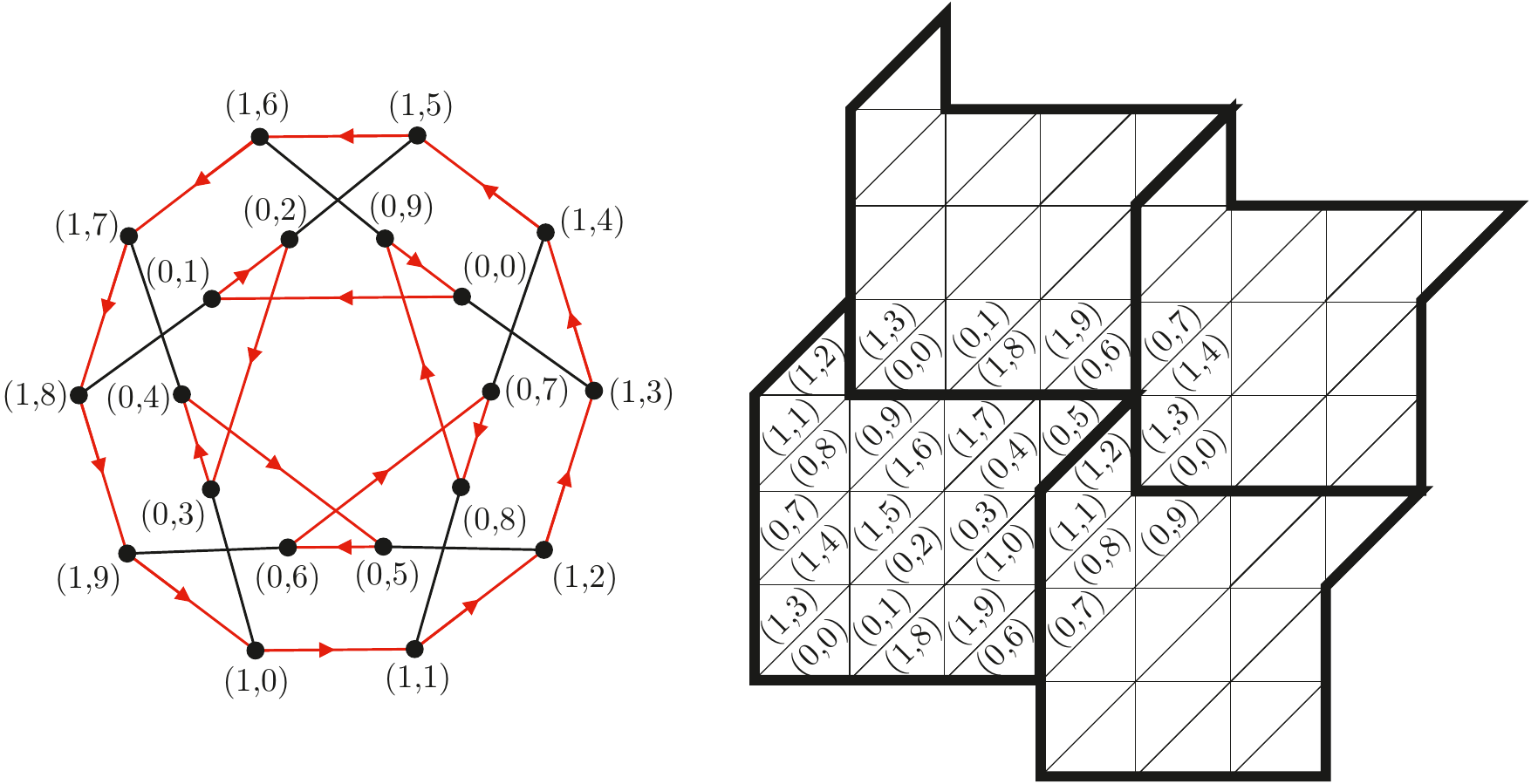}  
%\vskip -1.5cm
\caption{The chordal double ring mixed graph $CDRM(20,7)$, with diameter $k=6$, and its plane tessellation.}
\label{fig14}
\end{center}% \end{figure}
\end{figure}

\begin{remark}
As commented above, when $k$ is even but not a multiple of $4$,
the tessellation with the optimal tile of Figure \ref{fig11} does not give a directed ring with $n$ vertices. Instead, we obtain two directed rings on $m=n/2$ vertices each, which we call a {\em chordal double ring mixed graph} $CDRM(n,c)$. Then, we work on the group $\Z_2\times \Z_m$, and vertex $(\alpha,i)$ is adjacent to $(\alpha+1,i)$ (though an edge) and vertex $(\alpha,i+1)$ (through an arc). Thus, when $k$ is an even number of the form $k=4t+2$, the tile of Figure \ref{fig11} with $\ell=2t$ and order $n=2\ell^2+2$ corresponds to a 
$CRM(n,c)$ with $c=(\frac{k}{2}-1)^2+\frac{k}{2}$ (as in Theorem \ref{th:CRM}$(b)$). For instance, in the case $k=6$,
Figure \ref{fig14} shows the chordal ring mixed graph $CDRM(20,3)$ and its corresponding tessellation.
\end{remark}

%  \begin{figure}[!ht]
% 	\begin{center}
% 	\includegraphics[width=6cm]{DMC(20,3).jpg}  
% %\vskip -1.5cm
% 	\caption{corresponding to the mixed chordal double ring graph  $CDRM(20,3)$, with diameter $k=6$.}
% 		\label{fig15}
% 	\end{center}% \end{figure}
%  \end{figure}

In Table \ref{table6}, we show the values of the obtained chordal ring mixed graphs, given by Theorem \ref{th:CRM}, together with the upper bounds on the order of $CRM(n,c)$ and those from Theorem \ref{th:james}.
For example, notice that the last maximal $(1,1,5)$-mixed graph of Figure \ref{fig4} corresponds to the chordal ring mixed graph $CRM(18,5)$.

\begin{table}[ht]
\begin{center}
\begin{tabular}{|r||r|r|r|l|l|} \hline
$k$ & $n$ & $c$ & Th. \ref{th:CRM} & max order & Th. \ref{th:james}\\
& & &  & $CRM(n,c)$ & \\
\hline \hline
3 & 8 & 3 & $(a)$ &  8 & 8\\ \hline
4 & 10 & 3 & $(b)$ &  12 & 12\\ \hline
5 & 18 & 5 &  $(a)$ &  18 & 22\\ \hline
6 & 16 & 3 &  $(c1)$ &  24 & 36\\ \hline
7 & 32  &  7 &  $(a)$ &  32 &  60\\ \hline
8 & 34  & 13  & $(b)$ &  40 & 96\\ \hline
9 & 50  & 9  &  $(a)$ &   50 & 158\\ \hline
10 & 44 & 31 & $(c2)$ &  60 &  256\\ \hline
11 & 72 & 11  &  $(a)$ &  72 &  416\\ \hline
12 & 74 &  31 & $(b)$  &  84 & 674\\ \hline
13 & 98 & 13  &  $(a)$ &  98 &  1092 \\  \hline
14 & 88 & 19  & $(c1)$  &  112 & 1766\\ \hline
15 & 128 &  15 &  $(a)$ &   128 & 2860\\ \hline
16 & 139 &  57 & $(b)$  &  144 & 4628\\ \hline
17 & 162 & 17  &  $(a)$ &   162 & 7490\\ \hline
18 & 148 &  107 & $(c2)$  &  180 & 12120\\ \hline
19 & 200 & 19  &  $(a)$ &    200 & 19612 \\  \hline
20 & 202 & 91  & $(b)$  &  220 & 31732\\ \hline
21 & 242 &  21 &  $(a)$ &  242  & 51346\\ \hline
22 & 224 &  51 & $(c1)$  &  264 & 83080\\
\hline
\end{tabular}
\caption{Bounds for the chordal rings mixed graph $CRM(n,c)$. }
\label{table6}
\end{center}
\end{table}

\section{Upper bounds}\label{sec:upperbounds}

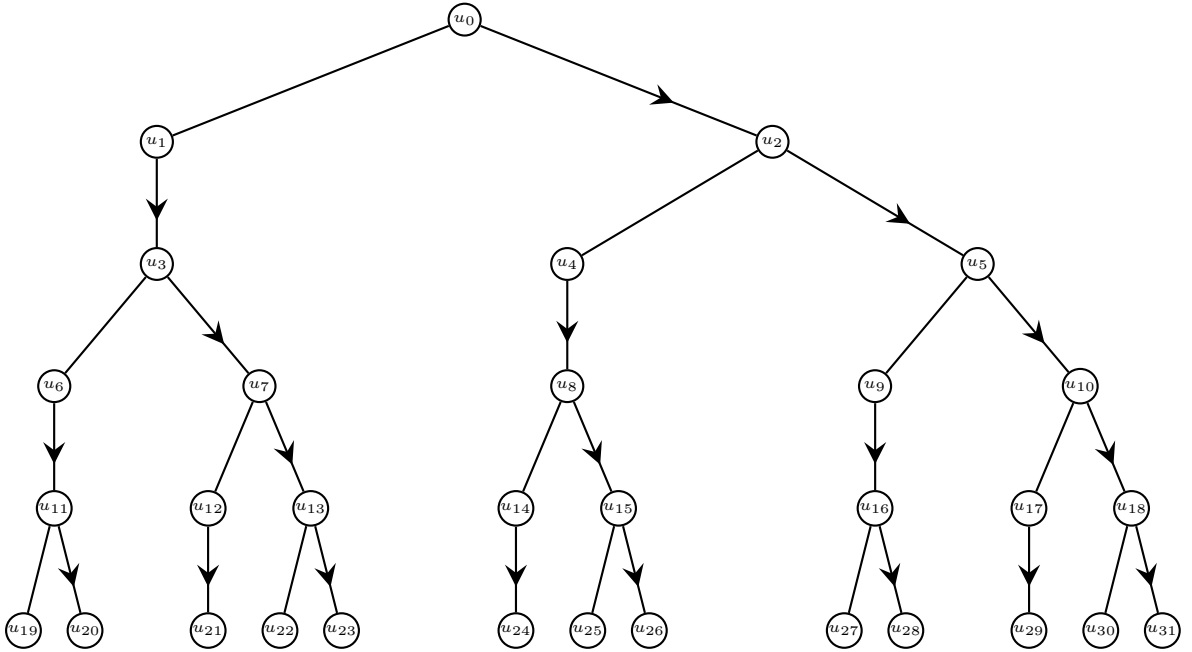
\begin{figure}\centering
\begin{tikzpicture}[middlearrow=stealth,x=0.2mm,y=-0.2mm,inner sep=0.1mm,scale=1.35,
thick,vertex/.style={circle,draw,minimum size=12,font=\tiny,fill=white},edge label/.style={fill=white}]
\tiny
\node at (0,0) [vertex] (u0) {$u_0$};
\node at (-150,60) [vertex] (u1) {$u_1$};
\node at (150,60) [vertex] (u2) {$u_2$};
\node at (-150,120) [vertex] (u3) {$u_3$};
\node at (50,120) [vertex] (u4) {$u_4$};
\node at (250,120) [vertex] (u5) {$u_5$};
\node at (-200,180) [vertex] (u6) {$u_6$};
\node at (-100,180) [vertex] (u7) {$u_7$};
\node at (50,180) [vertex] (u8) {$u_8$};
\node at (200,180) [vertex] (u9) {$u_9$};
\node at (300,180) [vertex] (u10) {$u_{10}$};
\node at (-200,240) [vertex] (u11) {$u_{11}$};
\node at (-125,240) [vertex] (u12) {$u_{12}$};
\node at (-75,240) [vertex] (u13) {$u_{13}$};
\node at (25,240) [vertex] (u14) {$u_{14}$};
\node at (75,240) [vertex] (u15) {$u_{15}$};
\node at (200,240) [vertex] (u16) {$u_{16}$};
\node at (275,240) [vertex] (u17) {$u_{17}$};
\node at (325,240) [vertex] (u18) {$u_{18}$};
\node at (-215,300) [vertex] (u19) {$u_{19}$};
\node at (-185,300) [vertex] (u20) {$u_{20}$};
\node at (-125,300) [vertex] (u21) {$u_{21}$};
\node at (-90,300) [vertex] (u22) {$u_{22}$};
\node at (-60,300) [vertex] (u23) {$u_{23}$};
\node at (25,300) [vertex] (u24) {$u_{24}$};
\node at (60,300) [vertex] (u25) {$u_{25}$};
\node at (90,300) [vertex] (u26) {$u_{26}$};
\node at (185,300) [vertex] (u27) {$u_{27}$};
\node at (215,300) [vertex] (u28) {$u_{28}$};
\node at (275,300) [vertex] (u29) {$u_{29}$};
\node at (310,300) [vertex] (u30) {$u_{30}$};
\node at (340,300) [vertex] (u31) {$u_{31}$};
\path
(u0) edge (u1)
(u0) edge [middlearrow] (u2)
(u1) edge [middlearrow] (u3)
(u2) edge (u4)
(u2) edge [middlearrow] (u5)
(u3) edge (u6)
(u3) edge [middlearrow] (u7)
(u4) edge [middlearrow] (u8)
(u5) edge (u9)
(u5) edge [middlearrow] (u10)
(u7) edge (u12)
(u8) edge (u14)
(u10) edge (u17)
(u6) edge [middlearrow] (u11)
(u7) edge [middlearrow] (u13)
(u8) edge [middlearrow] (u15)
(u9) edge [middlearrow] (u16)
(u10) edge [middlearrow] (u18)
(u11) edge (u19)
(u11) edge [middlearrow] (u20)
(u12) edge [middlearrow] (u21)
(u13) edge (u22)
(u13) edge [middlearrow] (u23)
(u14) edge [middlearrow] (u24)
(u15) edge (u25)
(u15) edge [middlearrow] (u26)
(u16) edge (u27)
(u16) edge [middlearrow] (u28)
(u17) edge [middlearrow] (u29)
(u18) edge (u30)
(u18) edge [middlearrow] (u31)
;
\end{tikzpicture}
\caption{The Moore tree for $r = z = 1$ and $k = 5$.}
\label{fig16}
\end{figure}

This section gives an upper bound that improves the Moore bound for bipartite mixed graphs. We concentrate on our case of $r=z=1$. The first step is to draw the Moore tree of depth $k$ for $r = z = 1$, starting at level 0. When a vertex at level $i$ has both an undirected edge and a directed arc to child vertices at level $i+1$, we draw the undirected neighbour to the left of the directed out-neighbour. See the Moore tree for diameter $k=5$ in Figure~\ref{fig16}.

We call a position in the Moore tree an \emph{arrow vertex} if it lies in level $i$ in the undirected branch of the tree for some $2 \leq i \leq k-1$ and arises as the endpoint of a directed arc from level $i-1$. Observe that if $u$ is an arrow vertex, then both elements of $N^-(u)$ occur in the undirected branch. As the directed out-neighbour of the root of the tree must reach $u$ by a mixed path of length at most $k$, at least one of these vertices of $N^-(u)$ must also occur in the directed branch of the tree, therefore contributing towards the defect of the graph. A complication is the fact that the in-neighbourhoods of arrow vertices can overlap. Hence, we require the smallest \emph{transversal} of the in-neighbourhoods of the vertices in the undirected branch. We, therefore, partition the set of the union of these in-neighbourhoods in a convenient way. If $v$ is a vertex at level $i$ %undirected branch 
for $1 \leq i \leq k-2$, then we denote by $v^{\rightarrow }$ the undirected neighbour of the directed out-neighbour of $v$, which appears at level $i+2$ of the undirected branch. 

By iterating this `$\rightarrow $' procedure, we obtain a chain of vertices that terminates at level $k-1$ or $k$. This chain is maximal if and only if it initiates at an arrow vertex or the undirected neighbour of the tree's root. These maximal chains are disjoint, and if there are $t$ vertices contained in the chain, then the smallest number of vertices from the chain that must also appear in the directed branch is given by $\left \lceil \frac{t}{3} \right \rceil $ (this is the domination number of the path $P_t$).

Now, we show how this argument can be applied to improve the bound for bipartite mixed graphs. We first focus on the case of even diameter $k = 2\kappa$. Consider maximal chains that occupy odd levels of the tree (for example, the chain beginning at the undirected neighbour of the root of the tree has vertices at levels $1,3,5,\dots,k-1$); we call this an \emph{odd level chain}. Recall that the Moore bound for bipartite mixed %parts 
graphs is twice the number of vertices on odd levels of the tree. By the preceding argument, the directed branch must contain the vertices from a transversal of the undirected branch, but by bipartiteness, each of the vertices of a transversal of an odd level chain must also lie on odd levels.

Tuite and Erskine \cite{te22} showed that there are $\eta _t = \frac{1}{2^{t-1}\sqrt{5}}((1+\sqrt{5})^{t-1}-(1-\sqrt{5})^{t-1})$ chains that start at level $t$ if $2\le t\le k-2$, and $\eta_1=1$. Counting over the odd level chains (and replacing the summation index $t$ by $2t-1$), we conclude that at least 
\begin{align*}
& \sum_{t=1}^{\kappa-1}\eta_{2t-1}\left\lceil\frac{\kappa-t+1}{3} \right\rceil \\
&=\left \lceil \frac{\kappa}{3} \right \rceil +\sum_ {t=2}^{\kappa-1}\frac{1}{2^{2t-2}\sqrt{5}}\Big((1+\sqrt{5})^{2t-2}-(1-\sqrt{5})^{2t-2}\Big)\left \lceil \frac{\kappa-t+1}{3} \right \rceil 
\end{align*}
vertices are repeated on the odd levels. As this applies to both partite sets, the defect is at least twice this figure.

By considering even level chains, we obtain the analogous conclusion for odd-diameter graphs.

\begin{theorem}
\label{th:james}
The order of a totally regular bipartite graph with undirected degree $r = 1$, directed degree $z = 1$, and diameter $k = 2\kappa$ is at most 
\[
M_B(1,1,2\kappa)- 2\left \lceil \frac{\kappa}{3} \right \rceil -2\sum _{t=2}^{\kappa-1}\frac{1}{2^{2t-2}\sqrt{5}}\Big((1+\sqrt{5})^{2t-2}-(1-\sqrt{5})^{2t-2}\Big)\left \lceil \frac{\kappa-t+1}{3} \right \rceil,
\]
for $\kappa\geq3$. For $\kappa=2$, the order is at most $M_B(1,1,2\kappa)- 2\left \lceil \frac{\kappa}{3} \right \rceil$.
If $k = 2\kappa+1$, then
\[ 
M_B(1,1,2\kappa+1)- 2\sum _{t=1}^{\kappa -1}\frac{1}{2^{2t-1}\sqrt{5}}\Big((1+\sqrt{5})^{2t-1}-(1-\sqrt{5})^{2t-1}\Big)\left \lceil \frac{\kappa-t+1}{3} \right \rceil,
\]
for $\kappa\geq2$.
\end{theorem}
For diameter $k = 4$, this gives an upper bound of 12, and for $k = 6$, an upper bound of 36. Moreover, for diameter $k=5$, our result gives the upper bound $22$. Except for this case, for which we know that the tight bound is, in fact, 18, this theorem gives all the other upper bounds in Table~\ref{table6}.

\section*{Acknowledgments}

We thank Nacho L\'opez for his valuable discussions about this subject.

\end{document}